\documentclass[10pt]{article}

\usepackage{a4wide}
\usepackage{amssymb}
\usepackage{amsfonts}
\usepackage{amsmath}
\input xy
\xyoption{arrow} \xyoption{matrix}

\usepackage{amsthm}
\usepackage{multirow}
\usepackage{marginnote}
\usepackage{a4wide}
\usepackage{amssymb}
\usepackage{amsfonts}
\usepackage{amsmath}
\usepackage{mathrsfs}
\usepackage{tikz}
\usepackage{tikz-cd}
\usetikzlibrary{arrows,matrix}
\usetikzlibrary{positioning}

\date{}

\newtheorem{proposition}{Proposition}[section]
\newtheorem{theorem}[proposition]{Theorem}
\newtheorem{lemma}[proposition]{Lemma}

\newtheorem{corollary}[proposition]{Corollary}

\def\Kdim{{\rm K.dim }\,}

\def\der{\partial }

\def\nFM0{{\nu }_{F,M_0}}
\def\nFN0{{\nu }_{F,N_0}}
\def\nGN0{{\nu }_{G,N_0}}

\def\N0{ {\bf N}_0 }

\def\g{\gamma}
\def\v{\varphi}
\def\ra{\rightarrow}

\def\Xpm{X^{\pm }}

\def\s{\sigma}
\def\Z{{\bf Z }}

\def\l1{{\lambda}_1}

\def\a{\alpha}
\def\a0{ {\alpha }_0}
\def\a1{ {\alpha }_1}

\def\l{\lambda}
\def\o{\omega}

\def\nFGM0{{\nu }_{F,G,M_0}}

\def\nFN0{{\nu}_{F,N_0}}

\def\sm{{\sigma}^m}

\def\sm1{{\sigma}^{-1}}

\def\smtp1{{\sigma}^{-t+1}}

\def\o{\omega }
\def\S1{S^{-1}}

\def\Xpm1{X^{\pm 1}_1}

\def\sPM1{{\sigma }^{\pm 1}}
\def\sMP1{{\sigma }^{\mp 1 }}

\def\d{\delta}

\def\di{{\rm d.ind}}

\def\L{\Lambda}

\def\G{\Gamma}

\def\CA{{\cal A}}

\def\CD{{\cal D}}

\def\Ytm1{Y^{t-1}}
\def\Yim1{Y^{i-1}}

\def\CK{{\cal K}}

\def\CN{{\cal N}}

\def\Aut{{\rm Aut}}

\def\ad{{\rm ad }}
\def\dim{{\rm dim }}

\def\ker{ {\rm ker } }

\def\SL2Z{ {\rm SL}_2({\bf Z}) }

\def\CR{ {\cal R}}

\def\th{ \theta }

\def\Gp1{ G^{1 , 1 } }
\def\P11{ P^{-1 , 1 } }
\def\Pp1{ P^{1 , 1 } }

\def\th{\theta}

\def\nCLsr{{}^\nu\kern-2pt {\cal L}^{\sigma , \rho  }}
\def\nP{{}^\nu \kern-2pt P}
\def\nL{{}^\nu\kern-2pt L}
\def\nLL{{}^\nu\kern-2pt \Lambda}
\def\nPsr{{}^\nu\kern-2pt P^{\sigma , \rho  }}
\def\nLsr{{}^\nu\kern-2pt L^{\sigma , \rho  }}
\def\nuCL{{}^\nu\kern-2pt  {\cal L}}
\def\nCLsr{{}^\nu\kern-2pt {\cal L}^{\sigma , \rho  }}
\def\nCL1m{{}^\nu\kern-2pt {\cal L}^{-1 , 1  }}
\def\x1nu{x^\frac{1}{\nu}}
\def\xm1nu{x^{-\frac{1}{\nu}}}

\def\CR{ {\cal R}}

\def\CN{{\cal N}}
\def\ra{\rightarrow }

\def\CB{{\cal B}}

\def\CI{{\cal I}}

\def\CC{ {\cal C}}

\def\nAM0{{\nu }_{{\cal A},M_0}}
\def\nAN0{{\nu }_{{\cal A},N_0}}

\def\Kdim{ {\rm Kdim } }

\def\End{ {\rm End }}

\def\CR{ {\cal R }}

\def\ad{ {\rm ad }}

\def\ga{\mathfrak{a}}

\def\gp{\mathfrak{p}}
\def\gq{\mathfrak{q}}

\def\Spec{{\rm Spec}}

\def\di!{\frac{\der^i}{i!}}
\def\dik!{\frac{\der^k_i}{k!}}

\def\CC{{\cal C}}

\def\Fp{\mathbb{F}_p}

\def\N{\mathbb{N}}

\def\0{\overline{0}}
\def\1{\overline{1}}

\def\Ln1{\L_{n,\overline{1}}}

\def\a1{a_{\overline{1}}}

\def\S{\Sigma}

\def\vn1{\overrightarrow{n-1}}

\def\Q{\mathbb{Q}}

\def\mA{\mathbb{A}}

\def\mS{\mathbb{S}}
\def\mJ{\mathbb{J}}
\def\mI{\mathbb{I}}

\def\cht{{\rm cht}}

\def\mF{\mathbb{F}}

\def\K1{{\rm K}_1}

\def\hmI1{\widehat{\mI_1}}
\def\tmI1{\widetilde{\mI_1}}
\def\tmJ1{\widetilde{\mJ_1}}
\def\hB1{\widehat{B_1}}
\def\hCB1{\widehat{\CB_1}}

\def\ga{\mathfrak{a}}

\def\Z{\mathbb{Z}}

\def\mB{\mathbb{B}}

\def\CW{{\cal W}}

\def\Es{{\rm Es}}
\def\Deg{{\rm Deg}}

\begin{document}

\author{V. V. \  Bavula and A. Al Khabyah 
}

\title{Bi-quadratic algebras on 3 generators with PBW: class II.1}

\maketitle
\begin{abstract}

In \cite{BigAlg-3gen},  an explicit description of bi-quadratic algebras  on three generators with PBW basis was obtained. There are four classes: I-IV. The aim of the paper is to study algebras that belong to one of the classes: class II.1. In this class, there are twelve distinct subclasses to consider. The aim of the paper is to study these algebras. In particular, the centre, the primes ideals, the automorphism groups and  the normal elements are explicitly described.  The Krull, classical Krull, global and Gelfand-Kirillov dimensions of the algebras and their tensor products are computed. \\



$${\bf Contents}$$
\begin{enumerate}
\item Introduction.
\item The centre of the algebra $A$.
\item The centre, the ideal structure  and the prime spectrum of the algebra $A_{x_1x_2}$.
\item The prime spectrum of the algebra $A$. 
\item The automorphism group of the algebra $A$. 

\item The Krull, global and Gelfand-Kirillov dimensions of the algebras $\CA = \bigotimes_{i=1}^n A_i$ where $A_i=A(q_i, \alpha_i ,\mu_i)$.
  
\end{enumerate}

\end{abstract}


\section{Introduction}
In this paper, $K$ is a field of arbitrary characteristic  char$(K)=p$ and $K^\times =K\backslash \{0\}$, $\Fp=\Z /p\Z$ is the field with $p$ elements,  modules means a left module, $q=q_1\in K\backslash \{ 0, 1\}$, $A=A(q, \alpha ,\mu )$ is the algebra as in Theorem \ref{A28Oct18}.(1) (if it is not stated otherwise), $\N =\{ 0,1, \ldots \}$ is the set of natural numbers and $\N_+=\{ 1, 2 \, \ldots \}$. \\

{\bf  Description  of bi-quadratic algebras on 3 generators with PBW basis,  an overview.}  \\ 

{\it Definition, \cite{BigAlg-3gen}}. Given $Q=(q_1,q_2,q_3)\in K^{\times 3},$ a $3\times 3$ matrix with entries in the field $K$,
\marginpar{AAabc}
\begin{equation}\label{AAabc}
\mA =
 \begin{pmatrix}
  a & b &  c\\
  \alpha & \beta & \g\\
  \lambda & \mu &  \nu\\
  
\end{pmatrix}
\end{equation}
and $\mB=(b_1,b_2,b_3)\in K^3$. A {\bf bi-quadratic algebra on 3 generators} $A=K[x_1,x_2,x_3;Q,\mA ,\mB]$ is an algebra  generated by the elements $x_1,x_2$ and $x_3$  subject to  the defining relations:\\
\marginpar{DRel1}
\begin{equation} \label{DRel1}
x_2x_1-q_1x_1x_2=ax_1+bx_2+cx_3+b_1,
\end{equation} 
\marginpar{DRel2}
\begin{equation}\label{DRel2}
x_3x_1-q_2x_1x_3=\alpha x_1+\beta x_2+\g x_3+b_2,
\end{equation}
\marginpar{DRel3}
\begin{equation}\label{DRel3}
x_3x_2-q_3x_2x_3=\lambda x_1+\mu x_2+\nu x_3+b_3.
\end{equation}
Clearly, $A=\sum_{\xi\in \N^3}Kx^\xi$ where $\xi = (\xi_1, \xi_2, \xi_3 )$ and $x^\xi = x_1^{\xi_1}x_2^{\xi_2}x_3^{\xi_3}$. We say that the algebra $A$ has {\bf PBW basis} if 
$A=\bigoplus_{\xi\in \N^3}Kx^\xi$.\\

In \cite{BigAlg-3gen},  an explicit description of bi-quadratic algebras  on 3 generators with PBW basis was obtained. The  algebras were partitioned into the following 4 classes:\\

I:  $q_1=q_2=q_3=1$, \\

II: $q_1\neq 1$, $q_2=q_3=1$,\\

III: $q_1\neq 1$, $q_2\neq 1$, $q_3=1$,\\ 

IV:  $q_1\neq 1$, $q_2\neq 1$, $q_3\neq 1$,\\


and algebras in each subclass were explicitly described. In \cite{BigAlg-3gen}, a classification (up to isomorphism) of algebras that belong to the first class ($q_1=q_2=q_3=1$) was given. Class II consists of two subclasses: II.1 (Theorem \ref{A28Oct18}.(1)) and II.2 (Theorem \ref{A28Oct18}.(2)). The aim of the paper is to study algebras (and their tensor products) that belong to class II.1. The class II.1 consists of twelve subclasses, see below.\\

{\bf  Description  of bi-quadratic algebras on 3 generators with PBW basis when $q_1\neq 1$, $q_2=q_3=1$.}

\begin{theorem}\label{A28Oct18}
\cite{BigAlg-3gen}  Suppose  that $q_1\neq 1$, $q_2=q_3=1$ and  $A=K[x_1,x_2,x_3;Q,\mA ,\mB]$. 
\begin{enumerate}
\item If $\mu +\alpha \neq 0$ then (up to $G_3$):
$$x_2x_1=q_1x_1x_2, \;\; x_3x_1= x_1(x_3+\alpha ) \;\; {\rm and}\;\; x_3x_2= x_2(x_3+\mu )$$
where $(\alpha , \mu ) \in \{ (1,\mu'), (0,1)\, | \, \mu'\in K\backslash \{ -1\} \}$. The elements $x_1$ and $x_2$ are regular normal elements of $A$. The algebra $A=K[x_3][x_2, x_1; \s , \tau , b=0, \rho =q_1]$ is a diskew polynomial algebra where $\tau (x_3) = x_3-\alpha$ and $ \s (x_3) = x_3-\mu$. Furthermore, the algebra $A= \CD [x_2, x_1; \s , \tau , a=h]$ is a GWA where $\CD = K[x_3][h; \tau \s]$ is a skew polynomial ring where $\tau \s (x_3)= x_3-\alpha - \mu$, $\s (h) = q_1 h$ and $\tau (h)= q_1^{-1}h$. The element $h$ is a regular normal element of $A$. 
\begin{enumerate}
\item The algebra $A=K[x_3][x_1;\s_1][x_2;\s_2]$ is an iterated skew polynomial algebra where $\s_1(x_3)=x_3-\alpha$, $\s_2 (x_3)=x_3-\mu$ and $\s_2(x_1)=q_1x_1$. 
\item The algebra homomorphism $A(q_1, \alpha =1, \mu =0)\ra A(q_1^{-1}, \alpha =0, \mu =1)$, $x_1\mapsto x_2$,  $x_2 \mapsto x_1$, $x_3 \mapsto x_3$ is an algebra isomorphism. 
\item In particular, the algebra homomorphism $A(q_1=-1, \alpha =1, \mu =0)\ra A(q_1=-1, \alpha =0, \mu =1)$, $x_1\mapsto x_2$,  $x_2 \mapsto x_1$, $x_3 \mapsto x_3$ is an algebra isomorphism. 
\end{enumerate}

\item If $\mu +\alpha = 0$ then (up to $G_3$):
$$x_2x_1=q_1x_1x_2+cx_3+b_1, \;\; x_3x_1= x_1(x_3+\alpha ) \;\; {\rm and}\;\; x_3x_2= x_2(x_3-\alpha )$$
where exactly one of the following cases occurs: 
\begin{enumerate}
\item $\alpha =0$ and  $c,b_1\in \{ 0,1 \}$, 
\item $\alpha =1$ and   $(c,b_1) \in \{ (0,0),  (1,1)\}$, 
\item  $\alpha =1$, $c=0$, $b_1=1$ or  $\alpha =1$, $c=1$,  $b_1\in K\backslash \{ 1\}$. 
\end{enumerate}
So, the cases (a)--(c) can be written as $(\alpha , c,  b_1)\in \{ 0,1\}^3$ or $\alpha =1$, $c=1$, $b_1\in K\backslash \{0,  1\}$.  The algebra $A=K[x_3][x_2, x_1; \s , \s^{-1} , b=cx_3+b_1, \rho =q_1]$ is a diskew polynomial algebra where $\s (x_3) = x_3+\alpha$. Furthermore, the algebra $A= K[x_3, h] [x_2, x_1; \s , a=h]$ is a classical GWA where $\s (x_3)= x_3+\alpha $ and $\s (h) = q_1h+cx_3+b_1 $.     
\end{enumerate}
\end{theorem}

Definitions of generalized Weyl algebras and diskew polynomial algebras are given in \cite{Bav-GWA-diskew}.\\

$\bullet$ {\em In this paper, $A=A(q, \alpha ,\mu )$ is an algebra that belongs to class II.1 (see Theorem \ref{A28Oct18}.(1)).} \\

Surprisingly, there are 12 different cases/subclasses 
 for the algebras $A$:
\begin{eqnarray*}
 (C1):& &\text{if  $q$ is not a root of unity and $p=0$},\\
(C2):&& \text{if  $q$ is not a root of unity, $p>0$ and} 
 \; (\alpha , \mu )\in \{ (1,\mu'), (0,1)\, | \, \mu'\in \mF_p\backslash \{ -1\}\},\\
(C3):&& \text{if  $q$ is not a root of unity, $p>0$, $\alpha =1$ and $\mu\not\in \mF_p $},\\
(C4):&& \text{if  $q$ is a primitive $n$'th  root of unity, $p=0$, $\alpha =1$ and $\mu =0$},\\
(C5):&& \text{if  $q$ is a primitive $n$'th  root of unity, $p=0$, $\alpha =0$ and $\mu =1$},\\
(C6):&& \text{if  $q$ is a primitive $n$'th  root of unity, $p=0$, $\alpha =1$ and $\mu =-\frac{\mu_1}{\mu_2}$ where $\mu_1,\mu_2\in \N_+$}\\
 &&  \text{such that $\mu_1\neq \mu_2$ and $(\mu_1,\mu_2)=1$},\\
(C7a):&& \text{if  $q$ is a primitive $n$'th  root of unity, $p=0$, $\alpha =1$ and $\mu =\frac{\mu_1}{\mu_2}$ where $\mu_1,\mu_2\in \N_+$}\\ 
&& \text{such that  $(\mu_1,\mu_2)=1$},\\
(C7b):&& \text{ if  $q$ is a primitive $n$'th  root of unity, $p=0$, $\alpha =1$ and $\mu\not\in \Q $},\\
(C8):&& \text{if  $q$ is  a primitive $n$'th  root of unity, $p>0$, $\alpha =1$, $\mu =0$ and $p\nmid n$}, \\
(C9): && \text{if  $q$ is  a primitive $n$'th  root of unity, $p>0$, $\alpha =0$, $\mu =1$ and $p\nmid n$,} \\
(C10):&& \text{if  $q$ is  a primitive $n$'th  root of unity, $p>0$, $\alpha =1$, $\mu\not\in \mF_p $  and $p\nmid n$},\\
(C11):&& \text{if  $q$ is  a primitive $n$'th  root of unity, $p>0$, $\alpha =1$, $\mu\in \mF_p \backslash \{ 0,-1\}=\{ 1,\ldots , p-2\}$  and $p\nmid n$.} 
\end{eqnarray*}
Algebras in these case  have different properties and  different techniques are used to study them. The paper is organized as follows. 

In Section \ref{CASETH5.1.(1)}, the centre of the algebra $A$ is described (Theorem \ref{C24Sep19}). It is proven that in the cases (C1)--(C10) the algebra $A$ is a free module over the centre and an explicit free basis is presented in each case (Theorem \ref{D27Aug20}). The algebra $A$ is a finitely generated module over the centre iff the algebra $A$ belongs to one of the cases C(8)--(C11) (Corollary \ref{a14Sep20}).

In Section \ref{SPECAX1X2}, properties of the localization $A_{x_1x_2}$ of the algebra $A$ at the powers of the regular normal element $x_1x_2$ are considered. They are used in the classification of the prime ideals of the algebra $A$. The centre of the algebra $A_{x_1x_2}$ is described (Theorem \ref{AC24Sep19}). Ideals  and prime ideals of the algebra $A_{x_1x_2}$ are classified  (Theorem \ref{15Sep20} and Theorem \ref{AA15Sep20}). The algebra $A_{x_1x_2}$ is a free module over its centre (Theorem \ref{ZD27Aug20}) and  an explicit set of free generators is presented (in each of 12 cases). The elements $x_1$ and $x_2$ are units of the algebra $A_{x_1x_2}$. Let $\o_{x_1}$ and $\o_{x_2}$ be the inner automorphism of $A_{x_1x_2}$ that are defined by the elements $x_1$ and $x_2$, respectively. Let $\ad_{x_3}$ be the inner derivation of $A_{x_1x_2}$ that is determined by the element $x_3$. In the triple $\o = (\o_{x_1},\o_{x_2}, \ad_{x_3})$, the $K$-linear maps  {\em commute}. The algebra $A_{x_1x_2}$ is a generalized weight $\o$-module (Proposition \ref{16Sep20}). The $\o$-eigenalgebra $\Es_{A_{x_1x_2}}(\o)$ is described (Corollary \ref{b16Sep20}). A simplicity criterion is given for the algebra $A_{x_1x_2}$ (Theorem \ref{A14Sep20}). The algebra $A_{x_1x_2}$ is a free left/right $\Es_{A_{x_1x_2}}(\o)$-module (Corollary \ref{b16Sep20}.(2)).

In Section \ref{SPECA},  a classification of  prime ideals of the algebra  $A$ is obtained (see (\ref{SpAAa}), (\ref{SpAAa3}) and   Theorem \ref{18Sep20}).

In Section \ref{AUTTH5.1.(1)}, an explicit  description  the automorphism group $\Aut_K(A)$ is given (Theorem \ref{A27Aug20}),   classifications of  $\o$-weight vectors  (Theorem \ref{20Sep20}) and normal elements  (Proposition \ref{A20Sep20}) of the algebra $A$ are presented. It is shown that every normal element of the algebra $A$ is an $\o$-weight vector, and vice versa (Proposition \ref{A20Sep20}).

In Section \ref{KRULLGLDIM},  the Krull, classical Krull,  global and Gelfand-Kirillov dimensions of the algebras  $\CA = \bigotimes_{i=1}^n A_i$ are computed where $A_i=A(q_i, \alpha_i ,\mu_i)$ (Proposition \ref{A23Aug20}.(3) and Theorem \ref{A24Sep19}).


\section{The centre of the algebra $A$}\label{CASETH5.1.(1)}

In this section,  the centre of the algebra $A$ is described (Theorem \ref{C24Sep19}). There are eleven distinct cases (C1)--(C11) (Theorem \ref{C24Sep19}). It is proven that in the cases (C1)--(C10) the algebra $A$ is a free module over the centre and an explicit free basis is presented in each case (Theorem \ref{D27Aug20}). The algebra $A$ is a finitely generated module over the centre iff the algebra $A$ belongs to one of the cases C(8)--(C11) (Corollary \ref{a14Sep20}).
\\

{\bf $\N^2$-grading of the algebra $A$.} By Theorem \ref{A28Oct18}.(1).(a), the algebra $A$
is an $\N^2$-graded algebra,
\begin{equation}\label{AN2gr}
A=\bigoplus_{\beta \in \N^2} K[x_3]x^\beta \;\; {\rm where}\;\; x^\beta = x_1^{\beta_1}x_2^{\beta_2}\;\; {\rm and}\;\; \beta =(\beta_1, \beta_2).
\end{equation}
The product in the algebra $A$ is given by the rule: For all elements $a,b\in K[x_3]$ and $\beta , \g \in \N^2$,
$$ ax^\beta \cdot bx^\g= a \o^\beta (b) q^{\beta_2\g_1} x^{\beta +\g}\;\; {\rm where}\;\; \o^\beta =\o_{x_1}^{\beta_1}\o_{x_2}^{\beta_2}.$$
So, the algebra $A$ is a skew crossed product of the (additive) monoid $\N^2$ with  the polynomial algebra $K[x_3]$. For all elements $\beta \in \N^2$, 
\begin{equation}\label{AN2gr1}
[x_3, x^\beta ] = (\beta_1\alpha +\beta_2\mu)x^\beta 
\end{equation}
where $[a,b]:=ab-ba$ is the commutator of elements $a$ and $b$. So, the direct sum (\ref{AN2gr}) is a direct sum of eigenspaces 
 for the inner derivation 
 $$\ad_{x_3}:A\ra A, \;\;a\mapsto \ad_{x_3}(a):=[x_3,a]$$ of the algebra $A$. \\

 {\bf The  algebra $K[x_3]^{\o_{x_1}, \o_{x_2}}$.} Given an algebra $R$, the set of all {\em regular elements} of $R$ (i.e. non-zero-divisors) is denoted by $\CC_R$. An element $x\in R$ is called a {\em normal element} of $R$ if $Rx=xR$. Let $\CN (R)$ be the set of normal elements of the algebra $R$. The set $\CN (R)$ is a multiplicative monoid that contains the group of units $R^\times $  and the centre  $Z(R)$ of the algebra $R$. Every regular normal element $x$ of $R$ determines an automorphism $\o_x$ of the algebra $R$ which is given by the rule:
\begin{equation}\label{xr=orx}
xr=\o_x(r)x\;\; {\rm for \; all}\;\; r\in R.
\end{equation}
If $x\in R^\times$ then $\o_x (r)=xrx^{-1}$ for all $r\in R$, i.e. $\o_x$ is an inner automorphism of $R$ that is determined by the element $x$. Let $V$ be a $K$-vector space. A $K$-linear map $\phi : V\ra V$ is called a {\em locally finite } (resp., {\em locally nilpotent}) map if for each element $v\in V$, $\dim_K(K[\phi ]v)<\infty$ (resp., $\phi^i v$ for all $i\gg 1$). A locally nilpotent map is a locally finite but not vice versa, in general. 

 The elements $x_1$ and $x_2$ of the algebra $A$ are regular normal elements. Clearly,
 the automorphisms $\o_{x_1}$ and $\o_{x_2}$ of the algebra $A$ act on its generators by the rule:
 \begin{eqnarray*}
\o_{x_1} &:& x_1\mapsto x_1, \;\;  x_2\mapsto q^{-1}x_2,\;\;  x_3\mapsto x_3-\alpha ,\\
\o_{x_2} &:& x_1\mapsto q x_1, \;\;  x_2\mapsto x_2,\;\;  x_3\mapsto x_3-\mu .
\end{eqnarray*}
In particular,  $\o_{x_1}\o_{x_2}= \o_{x_2} \o_{x_1}$. 
The automorphisms $\o_{x_1}$ and $\o_{x_2}$  of the algebra $A$ respect the $\N^2$-grading. Also the inner derivation $\ad_{x_3}$ of the algebra $A$ respects the $\N^2$-grading as follows from the equalities
$[x_3,x_1]=\alpha x_1$ and $[x_3,x_2]=\mu x_2$. {\em The $K$-linear maps $\o_{x_1}$, $\o_{x_2}$ and $\ad_{x_3}$ are commuting, locally finite maps that respect the $\N^2$-grading of the algebra $A$.}

Let us consider the polynomial 
\begin{equation}\label{fmu}
f_\mu (x_3):=\prod_{a,b\in \mF_p}(x_3-a-b\mu )=\prod_{b=0}^{p-1}\Big( x_3^p-x_3-b(\mu^p-\mu)\Big)= (x_3^p-x_3)^p-(\mu^p-\mu)^{p-1}(x_3^p-x_3).
\end{equation}
 In more detail, 
\begin{eqnarray*}
 f_\mu (x_3)&=&\prod_{b\in \mF_p}\prod_{a\in \mF_p}((x_3-b\mu )-a)=\prod_{b\in \mF_p}((x_3-b\mu )^p-(x_3-b\mu ))=\prod_{b\in \mF_p}\Big( x_3^p-x_3-b(\mu^p-\mu)\Big)\\
 &=& (x_3^p-x_3)^p-(\mu^p-\mu)^{p-1}(x_3^p-x_3).\\
\end{eqnarray*}

Let $G$ be the subgroup of the automorphism group $\Aut_K(K[x_3])$ which is generated by the restrictions of the automorphisms $\o_{x_1}$ and $\o_{x_2}$ to the subalgebra $K[x_3]$ of $A$. Let $K[x_3]^G=K[x_3]^{\o_{x_1}, \o_{x_2}}=\{ f\in K[x_3]\, | \, \o_{x_1}(f)=f, \o_{x_2}(f)=f\}$ be the algebra of $G$-invariants. 
 Proposition \ref{E24Sep19} is a description  of the algebra $K[x_3]^G$. 

  \begin{proposition}\label{E24Sep19}
$$K[x_3]^G=K[x_3]^{\o_{x_1}, \o_{x_2}}=\begin{cases}
K& \text{if $p=0$},\\
K[x_3^p-x_3]& \text{if $p>0$ and $(\alpha , \mu )\in \{ (1,\mu'), (0,1)\, | \, \mu'\in \mF_p\backslash \{ -1\}\} $},\\
K[f_\mu]& \text{if $p>0$, $\alpha =1$ and $\mu\not\in \mF_p $},
\end{cases}
$$
and  $f_\mu (x_3):=\prod_{a,b\in \mF_p}(x_3-a-b\mu )=\prod_{b=0}^{p-1}\Big( (x_3-b\mu )^p-(x_3-b\mu )\Big)=\prod_{a=0}^{p-1}\Big( (x_3-a )^p-\mu^{p-1}(x_3-a)\Big)$.
\end{proposition}

{\it Proof.} (i) {\em If $p=0$ then $K[x_3]^G=K$}: The condition that $\alpha +\mu \neq 0$ implies that either $\alpha \neq 0$ or $\mu \neq 0$ (or both). If $\alpha \neq 0$ (resp., $\mu \neq 0$) then $K[x_3]^{\o_{x_1}}=K$ (resp., $K[x_3]^{\o_{x_2}}=K$), and the statement (i)  follows since $K\subseteq K[x_3]^G\subseteq 
K[x_3]^{\o_{x_1}}=K$ (resp., $K\subseteq K[x_3]^G\subseteq 
K[x_3]^{\o_{x_2}}=K$).

(ii) {\em If $p>0$ and $(\alpha , \mu )\in \{ (1,\mu'), (0,1)\, | \, \mu'\in \mF_p\backslash \{ -1\}\} $ then $K[x_3]^G=K[x_3^p-x_3]$}: Since $\mF_p\alpha +\mF_p\mu = \mF_p$ and $\prod_{i=0}^{p-1} (x_3-i)=x_3^p-x_3\in K[x_3]^G$, the statement (ii) follows.

 (iii) {\em If $p>0$, $\alpha =1$ and $\mu\not\in \mF_p $ then $K[x_3]^G=K[f_\mu]$}: Since $f_\mu (x_3):=\prod_{a,b\in \mF_p}(x_3-a-b\mu )\in K[x_3]^G$ and the elements $\{ a+b\mu \, | \, a,b\in \mF_p\}$ are distinct elements of the field $K$, we have that $K[x_3]^G=K[f_\mu]$.  Indeed, let $g(x_3)\in K[x_3]^G$. We have to show that $g\in K[f_\mu]$. We use induction on the degree $n=\deg (g)$ of the polynomial $g$. The case $n=0$ is obvious. So, let $n>0$ and we assume that the result holds for all $G$-invariant  polynomials of degree $<n$. By replacing the polynomial $g(x_3)$ by $g(x_3)-g(0)$, we may assume that $g(x_3)=x_3h(x_3)$ for some polynomial $h(x_3)\in K[x_3]$. Now, for all numbers $a$ and $b$ such that $0\leq a,b\leq p-1$, 
 $$ g(x_3)=\o_{x_1}^a\o_{x_2}^b(g)=(x_3-a-b\mu )h(x_3-a-b\mu ).$$ 
 Hence, $g=\prod_{a,b\in \mF_p}(x_3-a-b\mu )
 \phi =f_\mu \phi$ for some polynomial $\phi=gf_\mu^{-1}\in K[x_3]^G$. Clearly, $\deg\, \phi = n-\deg\, f_\mu <n$. By induction $\phi \in K[f_\mu ]$, and so $g\in  K[f_\mu ]$, as required.

 Notice that $x_3^p-x_3=\prod_{a\in \mF_p}(x_3-a)$. Hence, $f_\mu = \prod_{b=0}^{p-1}\Big( (x_3-b\mu )^p-(x_3-b\mu )\Big)$.

 Notice that 
 $$\prod_{b\in \mF_p}(x_3-b\mu)=\mu^p \prod_{b\in \mF_p}(\mu^{-1}x_3-b)= \mu^p\Big( (\mu^{-1}x_3)^p-\mu^{-1}x_3\Big)=x_3^p-\mu^{p-1}x_3.
 $$
 
 Hence, $f_\mu=\prod_{a=0}^{p-1}\Big( (x_3-a )^p-\mu^{p-1}(x_3-a)\Big)$. $\Box$ \\

 {\bf The centralizer $C_A(x_1,x_2)$ of the elements $x_1$ and $x_2$ in the  algebra $A(q, \alpha ,\mu)$.} Let $S$ be a nonempty subset of an algebra $R$, the set $C_R(S):=\{ r\in R\, | \, rs=sr$ for all $s\in S \}$ is called the {\em centralizer} of $S$ in $R$. The centralizer $C_R(S)$ is a subalgebra of $R$.  Proposition \ref{D24Sep19} describes the centralizer $C_A(x_1,x_2)$.
 
 \begin{proposition}\label{D24Sep19}
The centralizer of the elements $x_1$ and $x_2$ is equal to 
\begin{enumerate}
\item $$C_A(x_1,x_2)=\begin{cases}
K[x_3]^G& \text{if $q$ is not a root of unity},\\
K[x_3]^G[x_1^n, x_2^n]& \text{if $q$ is a primitive $n$'th root of unity}.\\
\end{cases}
$$
\item $$C_A(x_1,x_2)=\begin{cases}
K& \text{if  $q$ is not a root of unity and $p=0$},\\
K[x_3^p-x_3]& \text{if  $q$ is not a root of unity, $p>0$ and} \\
 & (\alpha , \mu )\in \{ (1,\mu'), (0,1)\, | \, \mu'\in \mF_p\backslash \{ -1\}\} ,\\
K[f_\mu]& \text{if  $q$ is not a root of unity, $p>0$, $\alpha =1$ and $\mu\not\in \mF_p $},\\
K[x_1^n, x_2^n]& \text{if  $q$ is a primitive $n$'th  root of unity and $p=0$},\\
K[x_3^p-x_3][x_1^n, x_2^n]& \text{if  $q$ is  a primitive $n$'th  root of unity, $p>0$ and} \\
&  \text{$(\alpha , \mu )\in \{ (1,\mu'), (0,1)\, | \, \mu'\in \mF_p\backslash \{ -1\}\} $},\\
K[f_\mu][x_1^n, x_2^n]& \text{if  $q$ is  a primitive $n$'th  root of unity, $p>0$, $\alpha =1$ and $\mu\not\in \mF_p $},\\
\end{cases}
$$
\end{enumerate}
where the polynomial $f_\mu$ is defined in Proposition \ref{E24Sep19}. In particular, the centralizer $C_A(x_1,x_2)$ is a polynomial algebra in $m$ variables where $m=0,1,2,3$.
\end{proposition}

{\it Proof.} 1.  The elements $x_1$ and $x_2$ are regular and normal elements of the algebra $A$, and so  
$$C_A(x_1,x_2)=A^{\o_{x_1}, \o_{x_2}}.$$
The automorphisms $\o_{x_1}$ and $ \o_{x_2}$ of the algebra $A$ respect its $\N^2$-grading. Therefore, 
$$C_A(x_1,x_2)=\bigoplus_{\beta \in \N^2} \Big( K[x_3]x^\beta\Big)^{\o_{x_1}, \o_{x_2}},$$
and the statement follows.

2. Statement 2 follows from statement 1 and Proposition \ref{E24Sep19}.  $\Box$  \\

  {\bf The centre of algebra $A$.} For an algebra $R$ and its derivation $\d$, the set 
 $$R^\d :=\{ r\in R\, | \, \d (r)=0\}=\ker (\d )$$ is called the {\em algebra of $\d$-constants} in $A$. The elements $x_1$ and $x_2$ are regular and normal elements of the algebra $A$. Hence,  
\begin{equation}\label{ZAoder}
Z(A)=A^{\o_{x_1},\o_{x_2}, \ad_{x_3}} = C_A(x_1,x_2)^{\ad_{x_3}}=\bigoplus_{\beta \in \N^2} \bigg(\Big( K[x_3]x^\beta\Big)^{\o_{x_1}, \o_{x_2}}\bigg)^{\ad_{x_3}}.
\end{equation}
If $p>0$ and $q$ is a primitive $n$'th root of unity then $p\nmid n$ (since otherwise  $0=q^n-1=(q^\frac{n}{p}-1)^p$ implies $q^\frac{n}{p}=1$, a contradiction).  We use this fact in the proof of Theorem \ref{C24Sep19}. 
Theorem \ref{C24Sep19} is an explicit description of the centre of the algebra $A$. There are eleven distinct cases (C1)--(C11).

\begin{theorem}\label{C24Sep19}
Let $\S :=K[x_3^p-x_3][x_1^{pn}, x_1^{\xi_i n}x_2^{in},  x_2^{pn}]_{i=1, \ldots ,p-1}$ where $\xi_i$ is a unique natural number such that $0\leq \xi_i \leq p-1$ and $ \xi_i\equiv -i\mu \mod p$. Then
$$Z(A)=\begin{cases}
K& \text{(C1): if  $q$ is not a root of unity and $p=0$},\\
K[x_3^p-x_3]& \text{(C2): if  $q$ is not a root of unity, $p>0$ and} \\
 & (\alpha , \mu )\in \{ (1,\mu'), (0,1)\, | \, \mu'\in \mF_p\backslash \{ -1\}\} ,\\
K[f_\mu]& \text{(C3): if  $q$ is not a root of unity, $p>0$, $\alpha =1$ and $\mu\not\in \mF_p $},\\
K[x_2^n]& \text{(C4): if  $q$ is a primitive $n$'th  root of unity, $p=0$, $\alpha =1$ and $\mu =0$},\\
K[x_1^n]& \text{(C5): if  $q$ is a primitive $n$'th  root of unity, $p=0$, $\alpha =0$ and $\mu =1$},\\
K[x_1^{\mu_1n}x_2^{\mu_2n}]& \text{C(6): if  $q$ is a primitive $n$'th  root of unity, $p=0$, $\alpha =1$ and $\mu =-\frac{\mu_1}{\mu_2}$},\\
& \text{ where $\mu_1,\mu_2\in \N_+$ such that $\mu_1\neq \mu_2$ and $(\mu_1,\mu_2)=1$},\\
K& \text{(C7): if  $q$ is a primitive $n$'th  root of unity, $p=0$, $\alpha =1$ and $\mu\not\in \Q_-\cup \{ 0\}$},\\
K[x_3^p-x_3][x_1^{pn}, x_2^n]& \text{(C8): if  $q$ is  a primitive $n$'th  root of unity, $p>0$, $\alpha =1$, $\mu =0$ and $p\nmid n$}, \\
K[x_3^p-x_3][x_1^n, x_2^{pn}]& \text{(C9): if  $q$ is  a primitive $n$'th  root of unity, $p>0$, $\alpha =0$, $\mu =1$ and $p\nmid n$,} \\
K[f_\mu][x_1^{pn}, x_2^{pn}]& \text{(C10): if  $q$ is  a primitive $n$'th  root of unity, $p>0$, $\alpha =1$, $\mu\not\in \mF_p $  and $p\nmid n$},\\
\S& \text{(C11): if  $q$ is  a primitive $n$'th  root of unity, $p>0$, $\alpha =1$}, \\
& \text{$\mu\in \mF_p \backslash \{ 0,-1\}=\{ 1,\ldots , p-2\}$  and $p\nmid n$,} \\
\end{cases}
$$
where $\N_+=\N\backslash \{ 0\}$, $\Q_-=\{ q\in \Q\, | \, q<0\}$ and  $\xi_1=p-\mu $.
\end{theorem}

{\it Remarks.} 1. The case (C11) is vacuous if $p=2$ (since $\mF_2\backslash \{ 0,-1\}=\emptyset$). 

2. The algebra $\S$ is a commutative subalgebra of $A$ which is generated by $p+2$ elements that are defined in the theorem.

3. If $Z(A)\neq K$ then the centre $Z(A)$  uniquely determines the type of  the remaining nine cases in the Theorem \ref{C24Sep19}. 

4. The case (C7) consists of two subcases (C7a) and (C7b), see Theorem \ref{ZD24Sep19}. The ideal structure of the algebra $A$ is different in these two cases as in the case (C7a) the centre $Z(A_{x_1x_2})$ of the localized algebra $A_{x_1x_2}$ is not equal to $K$ and  ideals from the centre $Z(A_{x_1x_2})$
generate nontrivial ideals of the algebra $A$. \\

{\it Proof.} The theorem follows from (\ref{ZAoder}), Proposition \ref{D24Sep19}.(2) and (\ref{AN2gr1}).

 In more detail,  by Proposition \ref{D24Sep19}.(2) we have to consider 6 cases. The first three cases of   Proposition \ref{D24Sep19}.(2) are precisely  
the first three cases of the theorem (as in these cases $C_A(x_1,x_2)\subseteq Z(A)$, i.e. $C_A(x_1,x_2)= Z(A)$). 

The fourth case of Proposition \ref{D24Sep19}.(2) ($q$ is a primitive $n$'th  root of unity and $p=0$) is split into  
4 subcases, see cases 4-7 of the theorem. 
 The results follow from the equality $Z(A)=K[x_1^n,x_2^n]^{\ad_{x_3}}$ and (\ref{AN2gr1}).

For the fifth and sixth cases of   Proposition \ref{D24Sep19}.(2), we have that  $p\nmid n$ (since $q$ is a primitive $n$'th root of unity and $p>0$).


 The fifth case of Proposition \ref{D24Sep19}.(2) ($q$ is  a primitive $n$'th  root of unity, $p>0$ and
$(\alpha , \mu )\in \{ (1,\mu'), (0,1)\, | \, \mu'\in \mF_p\backslash \{ -1\}\} $) is split into  
3 subcases, see cases 8, 9 and 11 of the theorem. 
 The result follow from the equality $Z(A)=K[x_3^p-x_3][x_1^n,x_2^n]^{\ad_{x_3}}$ and (\ref{AN2gr1}).
 
 Finally,  the sixth case of Proposition \ref{D24Sep19}.(2) ($q$ is  a primitive $n$'th  root of unity, $p>0$, $\alpha =1$ and $\mu\not\in \mF_p $) is the 10'th  case of the theorem. 
  The result follow from the equality $Z(A)=K[f_\mu][x_1^n,x_2^n]^{\ad_{x_3}}$ and (\ref{AN2gr1}). $\Box$ \\

 {\bf The algebra $A$ is  free module over the centre in the cases  (C1)--(C10).}  Theorem \ref{D27Aug20} shows that the algebra $A$ is  free module over the centre in all cases of Theorem \ref{C24Sep19} but the one when $Z(A) = \S$ (i.e. $q$ is  a primitive $n$'th  root of unity, $p>0$, $\alpha =1$, $\mu\in \mF_p \backslash \{ 0,-1\}=\{ 1,\ldots , p-2\}$  and $p\nmid n$). 
 Explicit sets of free generators over $Z(A)$ are presented.

\begin{theorem}\label{D27Aug20}
 In the cases (C1)--(C10), the algebra $A=\bigoplus_{i\in I, \beta \in \CR} Z(A) x_3^ix^\beta$ is a free $Z(A)$-module 
 where the pair $(I,\CR )$ of sets $I=I_A\subseteq \N$ and $\CR = \CR_A\subseteq \N^2$ is given below where $\N_{<i}=\{ 0,1, \ldots , i-1\}$ and  $\mS (\mu_1n,\mu_2n) :=\N^2\backslash \Big((\mu_1n,\mu_2n)+\N^2\Big)$: 

 $$(I,\CR )=\begin{cases}
(\N, \N^2)& \text{(C1): if  $q$ is not a root of unity and $p=0$},\\
(\N_{<p},\N^2)& \text{(C2): if  $q$ is not a root of unity, $p>0$ and} \\
 & (\alpha , \mu )\in \{ (1,\mu'), (0,1)\, | \, \mu'\in \mF_p\backslash \{ -1\}\} ,\\
(\N_{<p^2},\N^2)& \text{(C3): if  $q$ is not a root of unity, $p>0$, $\alpha =1$ and $\mu\not\in \mF_p $},\\
(\N ,\N\times\N_{<n})& \text{(C4): if  $q$ is a primitive $n$'th  root of unity, $p=0$, $\alpha =1$ and $\mu =0$},\\
(\N ,\N_{<n}\times\N)& \text{(C5): if  $q$ is a primitive $n$'th  root of unity, $p=0$, $\alpha =0$ and $\mu =1$},\\
(\N ,\mS (\mu_1n, \mu_2n))& \text{(C6): if  $q$ is a primitive $n$'th  root of unity, $p=0$, $\alpha =1$ and $\mu =-\frac{\mu_1}{\mu_2}$},\\
& \text{ where $\mu_1,\mu_2\in \N_+$ such that $\mu_1\neq \mu_2$ and $(\mu_1,\mu_2)=1$},\\
(\N, \N^2)& \text{(C7): if  $q$ is a primitive $n$'th  root of unity, $p=0$, $\alpha =1$ and $\mu\not\in \Q_-\cup \{ 0\}$},\\

(\N_{<p} ,\N_{<pn}\times\N_{<n})& \text{(C8): if  $q$ is  a primitive $n$'th  root of unity, $p>0$, $\alpha =1$, $\mu =0$ and $p\nmid n$}, \\
(\N_{<p} ,\N_{<n}\times\N_{<pn})& \text{(C9): if  $q$ is  a primitive $n$'th  root of unity, $p>0$, $\alpha =0$, $\mu =1$ and $p\nmid n$,} \\
(\N_{<p^2} ,\N_{<pn}\times\N_{<pn})& \text{(C10): if  $q$ is  a primitive $n$'th  root of unity, $p>0$, $\alpha =1$, $\mu\not\in \mF_p $  and $p\nmid n$}.\\
\end{cases}
$$
In the case (C11), $A=\sum_{i\in I, \beta \in \CR} Z(A) x_3^ix^\beta =\bigoplus_{i\in I}x_3^i\Big( \sum_{ \beta \in \CR} Z(A) x^\beta\Big)$ where $I=\N_{<p}$ and $\CR =\N_{<pn}^2\backslash \Big( \bigcup_{i=0}^{p-1}((\xi_in, in)+\N_{<pn}^2)\Big)$.
\end{theorem}


{\it Proof.} The theorem follows at once from (\ref{AN2gr}) and Theorem \ref{C24Sep19}. $\Box$
\begin{corollary}\label{a14Sep20}
The algebra $A$ is a finitely generated module over its centre iff it belongs to the cases (C8)--(C11) of Theorem \ref{C24Sep19}.
\end{corollary}

{\it Proof}. By Theorem \ref{C24Sep19}, the algebra that belongs to the case (C11) case is finitely generated module over its centre. Now, the corollary follows from Theorem \ref{D27Aug20}. $\Box $\\

 {\bf The algebra $\Es_A(\o_{x_1},\o_{x_2})$.} Suppose that the algebra  $A$ belongs to the cases (C1)--(10). Then by Theorem \ref{D27Aug20}, $$A=\bigoplus_{i\in I_A, \beta \in \CR_A} Z(A) x_3^ix^\beta.$$ So, each  element $a\in A$ is a unique sum
 $a=\sum_{i\in I_A, \beta \in \CR_A} z_{i\beta} x_3^ix^\beta$ for some elements $z_{i\beta}\in Z(A)$. For $a\neq 0$, let 
 $${\rm Deg}_{x_3}(a):=\max \{ i \in I_A\, | \, z_{i\beta}\neq 0\;\; {\rm for \; some} \;\ \beta \in \CR_A\}. $$
   For all $z\in Z(A)$, $i\in I_A$ and $\beta \in \CR_A$,
\begin{equation}\label{oox1x}
\o_{x_1}(zx_3^ix^\beta )=q^{-\beta_2}z(x_3-\alpha )^ix^\beta\;\; {\rm and }\;\;  \o_{x_2}(zx_3^ix^\beta )=q^{\beta_1}z(x_3-\mu )^ix^\beta.
\end{equation}
    Clearly, for all elements $a\in A$, ${\rm Deg}_{x_3}(\o_{x_1}(a))={\rm Deg}_{x_3}(a)$ and ${\rm Deg}_{x_3}(\o_{x_2}(a))={\rm Deg}_{x_3}(a)$. 

Recall that $\o_{x_1}$ and $\o_{x_2}$ are commuting automorphisms of the algebra $A$.
 The $K$-linear span $ \Es_A(\o_{x_1},\o_{x_2})$  of all {\em common} eigenvectors of the automorphisms $\o_{x_1}$ and $\o_{x_2}$ is a subalgebra of $A$ which  is $\o_{x_1}-$ and $\o_{x_2}$-invariant.  The set $ \Es_A(\o_{x_1},\o_{x_2})$ is a {\em subalgebra} of $A$. We denote by the same symbols the restrictions of       the automorphisms $\o_{x_1}$ and $\o_{x_2}$ to the subalgebra $ \Es_A(\o_{x_1},\o_{x_2})$. 

 Corollary \ref{a28Aug20} is a description of the algebra $ \Es_A(\o_{x_1},\o_{x_2})$ where the algebra $A$ belongs to the cases (C1)--(C10). It also shows that the algebra $A$ is a free left/right $ \Es_A(\o_{x_1},\o_{x_2})$-module. 
 \begin{corollary}\label{a28Aug20}
 Suppose that the algebra $A$ belongs to one of the  cases (C1)--(C10) of Theorem \ref{C24Sep19}. Then 
\begin{enumerate}
\item $\Es_A(\o_{x_1},\o_{x_2})=\{a\in A \, | \,{\rm Deg}_{x_3}(a)\leq 0\}= \bigoplus_{\beta \in \CR_A}Z(A)x^\beta$, and  for all $z\in Z(A)$ and $\beta \in \CR_A$, $\o_{x_1}(zx^\beta )=q^{-\beta_2}zx^\beta$,  $\o_{x_2}(zx^\beta )=q^{\beta_1}zx^\beta$, and $\ad_{x_3}(zx^\beta)=(\alpha\beta_1+\mu \beta_2)zx^\beta$.
\item $A=\bigoplus_{i\in I_A}x_3^i\Es_A(\o_{x_1},\o_{x_2})=\bigoplus_{i\in I_A}\Es_A(\o_{x_1},\o_{x_2})x_3^i$ is a direct sum of $\Es_A(\o_{x_1},\o_{x_2})$-bimodules $x_3^i\Es_A(\o_{x_1},\o_{x_2})=\Es_A(\o_{x_1},\o_{x_2})x_3^i$. The algebra $A$ is a free left/right $ \Es_A(\o_{x_1},\o_{x_2})$-module. 
\end{enumerate}
\end{corollary}


{\it Proof}. 1. Let $E=\Es_A(\o_{x_1},\o_{x_2})$. The second equality of statement 1 is obvious. It remains to show that $E= \bigoplus_{\beta \in \CR_A}Z(A)x^\beta$. 
 This equality  follows from  Theorem \ref{D27Aug20} and (\ref{oox1x}). In more detail, by (\ref{oox1x}), the automorphisms $\o_{x_1}$ and $\o_{x_2}$ respect the $\N^2$-grading of the algebra $A$. Hence, by Theorem \ref{D27Aug20}, $E=\bigoplus_{\beta \in \CR}\Big(E\cap \bigoplus_{i\in I} Z(A)x^ix^\beta\Big)$.
By (\ref{oox1x}), for all elements $z\in Z(A)$ and $\beta \in \CR$, 
$$\o_{x_1}(zx^\beta )=q^{-\beta_2}zx^\beta\;\; {\rm  and }\;\;   \o_{x_2}(zx^\beta )=q^{\beta_1}zx^\beta .$$
If for some polynomial $\phi \in \bigoplus_{i\in I}Kx_3^i$,  $$ \o_{x_1}(\phi zx^\beta )=q^{-\beta_2}\phi zx^\beta \;\; {\rm and} \;\;  \o_{x_2}(\phi zx^\beta )=q^{\beta_1}\phi zx^\beta $$ then necessarily
 $\phi\in K[x_3]^G\subseteq Z(A)$, and the equality follows.
 
 2. Statement 2 follows from statement 1 and 
Theorem \ref{D27Aug20}. $\Box $

 Corollary \ref{a22Sep20} is a description of the algebra $ \Es_A(\o_{x_1},\o_{x_2})$ where the algebra $A$ belongs to the case (C11). It also shows that the algebra $A$ is a free left/right $ \Es_A(\o_{x_1},\o_{x_2})$-module. 
 \begin{corollary}\label{a22Sep20}
 Suppose that the algebra $A$ belongs to  the  case(C11)  of Theorem \ref{C24Sep19}. Let $\L' :=K[x_3^p-x_3][x_1^{pn},x_2^{pn}]$ and $I_A=\N_{<p}$. Then $\L'\subseteq Z(A)$ and 
\begin{enumerate}
\item $\Es_A(\o_{x_1},\o_{x_2})= \bigoplus_{\beta \in \N_{<pn}^2}\L'x^\beta$, and  for all $z\in \L'$ and $\beta \in \N^2_{<pn}$, $\o_{x_1}(zx^\beta )=q^{-\beta_2}zx^\beta$,  $\o_{x_2}(zx^\beta )=q^{\beta_1}zx^\beta$, and $\ad_{x_3}(zx^\beta)=(\alpha\beta_1+\mu \beta_2)zx^\beta$.
\item $A=\bigoplus_{i\in I_A}x_3^i\Es_A(\o_{x_1},\o_{x_2})=\bigoplus_{i\in I_A}\Es_A(\o_{x_1},\o_{x_2})x_3^i$ is a direct sum of $\Es_A(\o_{x_1},\o_{x_2})$-bimodules $x_3^i\Es_A(\o_{x_1},\o_{x_2})=\Es_A(\o_{x_1},\o_{x_2})x_3^i$. The algebra $A$ is a free left/right $ \Es_A(\o_{x_1},\o_{x_2})$-module. 
\end{enumerate}
\end{corollary}


{\it Proof}. 1. The inclusion $\L'\subseteq Z(A)$ is obvious. Let $E=\Es_A(\o_{x_1},\o_{x_2})$. Notice that $A= \bigoplus_{i\in I
_A, \beta \in \N_{<pn}^2}\L'x_3^ix^\beta$ and each direct summand is $\o_{x_i}$-invariant where $i=1,2$. Therefore,  $E=\bigoplus_{i\in I_A, \beta \in \N_{<pn}^2}\Big(E\cap\L'x_3^ix^\beta \Big)$.
By (\ref{oox1x}), for all elements $z\in \L'$ and $\beta \in \CR$, 
$$\o_{x_1}(zx^\beta )=q^{-\beta_2}zx^\beta\;\; {\rm  and }\;\;   \o_{x_2}(zx^\beta )=q^{\beta_1}zx^\beta .$$
If for some polynomial $\phi \in \bigoplus_{i\in I_A}Kx_3^i$,  $$ \o_{x_1}(\phi zx^\beta )=q^{-\beta_2}\phi zx^\beta \;\; {\rm and} \;\;  \o_{x_2}(\phi zx^\beta )=q^{\beta_1}\phi zx^\beta $$ then necessarily
 $\phi\in K[x_3]^G\subseteq \L'$, and the equality follows.
 
 2. Statement 2 follows from statement 1 and 
the equality $A= \bigoplus_{i\in I
_A, \beta \in \N_{<pn}^2}\L'x_3^ix^\beta$. $\Box $


\section{The centre, the ideal structure  and the prime spectrum of the algebra $A_{x_1x_2}$}\label{SPECAX1X2}

The aim of this section is for the algebra $A_{x_1x_2}$ to describe its  centre (Theorem \ref{AC24Sep19}) and to classify its ideals  (Theorem \ref{15Sep20}) and prime ideals  (Theorem \ref{AA15Sep20}). We show that the algebra $A_{x_1x_2}$ is a free module over its centre (Theorem \ref{ZD27Aug20}) and present an explicit set of free generators (in each of 12 cases). The algebra $A_{x_1x_2}$ is a generalized weight $\o$-module (Proposition \ref{16Sep20}) where $\o =(\o_{x_1},\o_{x_2}, \ad_{x_3})$. The $\o$-eigenalgebra $\Es_{A_{x_1x_2}}(\o)$ is described (Corollary \ref{b16Sep20}). A simplicity criterion is given for the algebra $A_{x_1x_2}$ (Theorem \ref{A14Sep20}). The algebra $A_{x_1x_2}$ is a free left/right $\Es_{A_{x_1x_2}}(\o)$-module (Corollary \ref{b16Sep20}.(2)).  \\

{\bf The algebra $A_{x_1x_2}$ and its $\Z^2$-grading.} Recall that the elements $x_1$ and $x_2$ are regular normal elements of the algebra $A=A(q, \alpha , \mu)_{x_1x_2}$. Hence, so is their product $x_1x_2$. Let $A_{x_1x_2}$ be the localization of the algebra $A$ at the powers of the element $x_1x_2$. The set S$=\{x_1^ix_2^j\, | \, i,j\geq 0\}$ is an Ore set of the domain $A$. Clearly, $A_{x_1x_2}\simeq S^{-1}A\simeq AS^{-1}$.

The algebra $A_{x_1x_2}$
is an $\Z^2$-graded algebra,
\begin{equation}\label{ZAN2gr} 
A_{x_1x_2}=\bigoplus_{\beta \in \Z^2} K[x_3]x^\beta \;\; {\rm where}\;\; x^\beta = x_1^{\beta_1}x_2^{\beta_2}\;\; {\rm and}\;\; \beta =(\beta_1, \beta_2).
\end{equation}
The product in the algebra $A$ is given by the rule: For all elements $a,b\in K[x_3]$ and $\beta , \g \in \Z^2$,
$$ ax^\beta \cdot bx^\g= a \o^\beta (b) q^{\beta_2\g_1} x^{\beta +\g}\;\; {\rm where}\;\; \o^\beta =\o_{x_1}^{\beta_1}\o_{x_2}^{\beta_2}.$$
So, the algebra $A_{x_1x_2}$ is a skew crossed product of the group $\Z^2$ with  the polynomial algebra $K[x_3]$. For all elements $\beta \in \Z^2$, 
\begin{equation}\label{ZAN2gr1}
[x_3, x^\beta ] = (\beta_1\alpha +\beta_2\mu)x^\beta .
\end{equation}
So, the direct sum (\ref{ZAN2gr}) is a direct sum of eigenspaces 
 for the inner derivation 
 $\ad_{x_3}$ of the algebra $A_{x_1x_2}$.  
 The elements $x_1$ and $x_2$ of the algebra $A_{x_1x_2}$ are units. Clearly,
 the automorphisms $\o_{x_1}$ and $\o_{x_2}$ of the algebra $A_{x_1x_2}$ are the inner automorphisms of the algebra $A_{x_1x_2}$ that are determined by the units $x_1$ and $x_2$. {\em The $K$-linear maps $\o_{x_1}$, $\o_{x_2}$ and $\ad_{x_3}$ are commuting, locally finite maps that respect the $\Z^2$-grading of the algebra $A_{x_1x_2}$.}\\

{\bf The centralizer $C_{A_{x_1x_2}}(x_1,x_2)$ of the elements $x_1$ and $x_2$ in the  algebra $A_{x_1x_2}$.}   Proposition \ref{ZD24Sep19} describes the centralizer $C_A(x_1,x_2)$.
 
 \begin{proposition}\label{ZD24Sep19}
The centralizer of the elements $x_1$ and $x_2$ is equal to 
\begin{enumerate}
\item $$C_{A_{x_1x_2}}(x_1,x_2)=\begin{cases}
K[x_3]^G& \text{if $q$ is not a root of unity},\\
K[x_3]^G[x_1^{\pm n}, x_2^{\pm n}]& \text{if $q$ is a primitive $n$'th root of unity}.\\
\end{cases}
$$
\item $$C_{A_{x_1x_2}}(x_1,x_2)=\begin{cases}
K& \text{if  $q$ is not a root of unity and $p=0$},\\
K[x_3^p-x_3]& \text{if  $q$ is not a root of unity, $p>0$ and} \\
 & (\alpha , \mu )\in \{ (1,\mu'), (0,1)\, | \, \mu'\in \mF_p\backslash \{ -1\}\} ,\\
K[f_\mu]& \text{if  $q$ is not a root of unity, $p>0$, $\alpha =1$ and $\mu\not\in \mF_p $},\\
K[x_1^{\pm n}, x_2^{\pm n}]& \text{if  $q$ is a primitive $n$'th  root of unity and $p=0$},\\
K[x_3^p-x_3][x_1^{\pm n}, x_2^{\pm n}]&  \text{if  $q$ is  a primitive $n$'th  root of unity, $p>0$ and} \\
&  \text{$(\alpha , \mu )\in \{ (1,\mu'), (0,1)\, | \, \mu'\in \mF_p\backslash \{ -1\}\} $},\\
K[f_\mu][x_1^{\pm n}, x_2^{\pm n}]& \text{if  $q$ is  a primitive $n$'th  root of unity, $p>0$, $\alpha =1$ and $\mu\not\in \mF_p $},\\
\end{cases}
$$
\end{enumerate}
where the polynomial $f_\mu$ is defined in Proposition \ref{E24Sep19}. 
\end{proposition}

{\it Proof.} 1.  Recall that the commuting automorphisms  $\o_{x_1}$ and $ \o_{x_2}$ of the algebra $A_{x_1x_2}$ respect its  $\Z^2$-grading. Therefore, 
$$C_{A_{x_1x_2}}(x_1,x_2)=A_{x_1x_2}^{\o_{x_1}, \o_{x_2}}=\bigoplus_{\beta \in \Z^2} \Big( K[x_3]x^\beta\Big)^{\o_{x_1}, \o_{x_2}},$$
and the statement follows.

2. Statement 2 follows from statement 1 and Proposition \ref{E24Sep19}.  $\Box$  \\

  {\bf The centre of algebra $A_{x_1x_2}$.} Theorem \ref{AC24Sep19} is an explicit description of the centre of the algebra $A_{x_1x_2}$.  Notice that 
\begin{equation}\label{ZAoder1}
Z(A_{x_1x_2})=A_{x_1x_2}^{\o_{x_1}, \o_{x_2},\ad_{x_3}} = C_{A_{x_1x_2}}(x_1,x_2)^{\ad_{x_3}}=\bigoplus_{\beta \in \Z^2} \bigg(\Big( K[x_3]x^\beta\Big)^{\o_{x_1}, \o_{x_2}}\bigg)^{\ad_{x_3}}.
 \end{equation}
 \begin{theorem}\label{AC24Sep19}
 Let $\S':=K[x_3^p-x_3][x_1^{\pm pn}, (x_1^{-\mu n}x_2^n)^{\pm 1}, x_2^{\pm pn}]$. Then 
$$Z(A_{x_1x_2})=\begin{cases}
K& \text{(C1): if  $q$ is not a root of unity and $p=0$},\\
K[x_3^p-x_3]& \text{(C2): if  $q$ is not a root of unity, $p>0$ and} \\
 & (\alpha , \mu )\in \{ (1,\mu'), (0,1)\, | \, \mu'\in \mF_p\backslash \{ -1\}\} ,\\
K[f_\mu]& \text{(C3): if  $q$ is not a root of unity, $p>0$, $\alpha =1$ and $\mu\not\in \mF_p $},\\
K[x_2^{\pm n}]& \text{(C4): if  $q$ is a primitive $n$'th  root of unity, $p=0$, $\alpha =1$ and $\mu =0$},\\
K[x_1^{\pm n}]& \text{(C5): if  $q$ is a primitive $n$'th  root of unity, $p=0$, $\alpha =0$ and $\mu =1$},\\
K[(x_1^{\mu_1n}x_2^{\mu_2n})^{\pm 1}]& \text{(C6): if  $q$ is a primitive $n$'th  root of unity, $p=0$, $\alpha =1$ and $\mu =-\frac{\mu_1}{\mu_2}$},\\
& \text{ where $\mu_1,\mu_2\in \N_+$ such that $\mu_1\neq \mu_2$ and $(\mu_1,\mu_2)=1$},\\
K[(x_1^{\mu_1n}x_2^{-\mu_2n})^{\pm 1}]& \text{(C7a): if  $q$ is a primitive $n$'th  root of unity, $p=0$, $\alpha =1$ and $\mu =\frac{\mu_1}{\mu_2}$},\\
& \text{ where $\mu_1,\mu_2\in \N_+$ such that  $(\mu_1,\mu_2)=1$},\\
K& \text{(C7b): if  $q$ is a primitive $n$'th  root of unity, $p=0$, $\alpha =1$ and $\mu\not\in \Q $},\\
K[x_3^p-x_3][x_1^{\pm pn}, x_2^{\pm n}]& \text{(C8): if  $q$ is  a primitive $n$'th  root of unity, $p>0$, $\alpha =1$, $\mu =0$ and $p\nmid n$}, \\
K[x_3^p-x_3][x_1^{\pm n}, x_2^{\pm pn}]& \text{(C9): if  $q$ is  a primitive $n$'th  root of unity, $p>0$, $\alpha =0$, $\mu =1$ and $p\nmid n$,} \\
K[f_\mu][x_1^{\pm pn}, x_2^{\pm pn}]& \text{(C10): if  $q$ is  a primitive $n$'th  root of unity, $p>0$, $\alpha =1$, $\mu\not\in \mF_p $  and $p\nmid n$},\\
\S'& \text{(C11): if  $q$ is  a primitive $n$'th  root of unity, $p>0$, $\alpha =1$}, \\
& \text{$\mu\in \mF_p \backslash \{ 0,-1\}=\{ 1,\ldots , p-2\}$  and $p\nmid n$.} \\
\end{cases}
$$
\end{theorem}

{\it Remarks.} 1. The case (C7) is split into two sub-cases (C7a) and (C7b).

2. The (C11) is vacuous if $p=2$ (since $\mF_2\backslash \{ 0,-1\}=\emptyset$). The case (C7b) is vacuous if $K=\Q$. 

3. If $Z(A_{x_1x_2})\neq K$ then the centre $Z(A_{x_1x_2})$  uniquely determines the type of  the remaining nine cases in the Theorem \ref{AC24Sep19}. \\

{\it Proof.} The theorem follows from (\ref{ZAoder1}), Proposition \ref{ZD24Sep19}.(2) and (\ref{ZAN2gr1}).

 In more detail,  by Proposition \ref{ZD24Sep19}.(2) we have to consider 6 cases. The first three cases of   Proposition \ref{ZD24Sep19}.(2) are precisely  
the first three cases of the theorem (as in these cases $C_{A_{x_1x_2}}(x_1,x_2)\subseteq Z(A_{x_1x_2})$, i.e. $C_{A_{x_1x_2}}(x_1,x_2)= Z(A_{x_1x_2})$). 

The fourth case of Proposition \ref{ZD24Sep19}.(2) ($q$ is a primitive $n$'th  root of unity and $p=0$) is split into  
5 subcases, see the cases (C4)--(C7a) and (C7b) of the theorem. 
 The results follow from the equality $Z(A_{x_1x_2})=K[x_1^{\pm n},x_2^{\pm n}]^{\ad_{x_3}}$ and (\ref{ZAN2gr1}).

For the fifth and sixth cases of   Proposition \ref{ZD24Sep19}.(2), we have that  $p\nmid n$ (since $q$ is a primitive $n$'th root of unity and $p>0$).


 The fifth case of Proposition \ref{ZD24Sep19}.(2) ($q$ is  a primitive $n$'th  root of unity, $p>0$ and
$(\alpha , \mu )\in \{ (1,\mu'), (0,1)\, | \, \mu'\in \mF_p\backslash \{ -1\}\} $) is split into  
3 subcases, see cases (C8), (C9) and (C11) of the theorem. 
 The result follow from the equality $Z(A_{x_1x_2})=K[x_3^p-x_3][x_1^{\pm n},x_2^{\pm n}]^{\ad_{x_3}}$ and (\ref{ZAN2gr1}).
 
 Finally,  the sixth case of Proposition \ref{ZD24Sep19}.(2) ($q$ is  a primitive $n$'th  root of unity, $p>0$, $\alpha =1$ and $\mu\not\in \mF_p $) is the case (C10) of the theorem. 
  The result follow from the equality $Z(A_{x_1x_2})=K[f_\mu][x_1^{\pm n},x_2^{\pm n}]^{\ad_{x_3}}$ and (\ref{ZAN2gr1}). $\Box$ 
  
The following obvious lemma is used in several proofs  in order to reduce the number of cases to deal with.  In particular, by Lemma \ref{a16Sep20}.(1) (resp., Lemma \ref{a16Sep20}.(1)), the cases (C4) and (C5) (resp., (C6) and (C7a)) are the same up to replacing $q$ by $q^{-1}$. Similarly, by Lemma \ref{a16Sep20}.(1), the cases (C8) and (C9)  are the same up to replacing $q$ by $q^{-1}$.
  
\begin{lemma}\label{a16Sep20}
\begin{enumerate}
\item The  homomorphism $A_{x_1x_2}(q, \alpha , \mu ) \ra A_{x_1x_2}(q^{-1}, \mu ,  \alpha  )$, $x_1\mapsto x_2$, $x_2\mapsto x_1$ is an  isomorphism. 
\item The  homomorphism $A_{x_1x_2}(q, \alpha , \mu ) \ra A_{x_1x_2}(q^{-1}, -\alpha ,  \mu  )$, $x_1\mapsto x_1^{-1}$, $x_2\mapsto x_2$ is an  isomorphism.
\end{enumerate}
\end{lemma}

{\it Proof}. Straightforward. $\Box $ \\

{\bf The algebra $A_{x_1x_2}$ is  free module over the centre.}  Theorem \ref{ZD27Aug20} shows that the algebra $A_{x_1x_2}$ is  free module over its centre and  
 explicit sets of free generators are presented.

\begin{theorem}\label{ZD27Aug20}
The algebra $A_{x_1x_2}=\bigoplus_{i\in I, \beta \in \CR} Z(A_{x_1x_2}) x_3^ix^\beta$ is a free $Z(A_{x_1x_2})$-module 
 where the pair $(I,\CR )$ of sets $I=I_{A_{x_1x_2}}\subseteq \N$ and $\CR = \CR_{A_{x_1x_2}}\subseteq \Z^2$ is given below where $\N_{<i}=\{ 0,1, \ldots , i-1\}$, $\mS' (\mu_1n,\mu_2n) :=\N\times \Z\backslash \Big((\mu_1n,\mu_2n)+\N\times \Z\Big)$ and  $\mS'' (-\mu_1n,\mu_2n) :=\Z\times \N\backslash \Big((-\mu_1n,\mu_2n)+\Z\times \N\Big)$: 

 $$(I,\CR )=\begin{cases}
(\N, \Z^2)& \text{(C1): if  $q$ is not a root of unity and $p=0$},\\
(\N_{<p},\Z^2)& \text{(C2): if  $q$ is not a root of unity, $p>0$ and} \\
 & (\alpha , \mu )\in \{ (1,\mu'), (0,1)\, | \, \mu'\in \mF_p\backslash \{ -1\}\} ,\\
(\N_{<p^2},\Z^2)& \text{(C3): if  $q$ is not a root of unity, $p>0$, $\alpha =1$ and $\mu\not\in \mF_p $},\\
(\N ,\Z\times\N_{<n})& \text{(C4): if  $q$ is a primitive $n$'th  root of unity, $p=0$, $\alpha =1$ and $\mu =0$},\\
(\N ,\N_{<n}\times\Z)& \text{(C5): if  $q$ is a primitive $n$'th  root of unity, $p=0$, $\alpha =0$ and $\mu =1$},\\
(\N ,\mS' (\mu_1n, \mu_2n))& \text{(C6): if  $q$ is a primitive $n$'th  root of unity, $p=0$, $\alpha =1$ and $\mu =-\frac{\mu_1}{\mu_2}$},\\
& \text{ where $\mu_1,\mu_2\in \N_+$ such that $\mu_1\neq \mu_2$ and $(\mu_1,\mu_2)=1$},\\
(\N ,\mS'' (-\mu_1n, \mu_2n)) &\text{(C7a): if  $q$ is a primitive $n$'th  root of unity, $p=0$, $\alpha =1$ and $\mu =\frac{\mu_1}{\mu_2}$},\\
& \text{ where $\mu_1,\mu_2\in \N_+$ such that $(\mu_1,\mu_2)=1$},\\
(\N , \Z^2)& \text{(C7b): if  $q$ is a primitive $n$'th  root of unity, $p=0$, $\alpha =1$ and $\mu\not\in \Q$},\\
(\N_{<p} ,\N_{<pn}\times\N_{<n})& \text{(C8): if  $q$ is  a primitive $n$'th  root of unity, $p>0$, $\alpha =1$, $\mu =0$ and $p\nmid n$}, \\
(\N_{<p} ,\N_{<n}\times\N_{<pn})& \text{C(9): if  $q$ is  a primitive $n$'th  root of unity, $p>0$, $\alpha =0$, $\mu =1$ and $p\nmid n$,} \\
(\N_{<p^2} ,\N_{<pn}\times\N_{<pn})& \text{(C10): if  $q$ is  a primitive $n$'th  root of unity, $p>0$, $\alpha =1$, $\mu\not\in \mF_p $  and $p\nmid n$},\\
(\N_{<p} ,\N_{<pn}\times\N_{<n})& \text{(C11): if  $q$ is  a primitive $n$'th  root of unity, $p>0$, $\alpha =1$}, \\
& \text{$\mu\in \mF_p \backslash \{ 0,-1\}=\{ 1,\ldots , p-2\}$  and $p\nmid n$.}
\end{cases}
$$
The algebra $A$ is a finitely generated module over its centre iff it belongs to  the last four cases. 
\end{theorem}


{\it Proof.} The theorem follows at once from (\ref{ZAN2gr}) and Theorem \ref{AC24Sep19}. $\Box$\\

Let $V$ be  a $K$-vector space and $a_1, \ldots , a_m\in \End_K(V)$ be commuting $K$-linear maps.  A nonzero vector $v\in V$ is called a {\em common eigenvector} of the maps $a_1, \ldots , a_m$ if $a_iv=\l_iv$ for $i=1, \ldots , m$ and some scalars $\l_i\in K$. 
 The $m$-tuple $\l = (\l_1, \ldots , \l_m)$ is called a {\em weight} for $a=(a_1, \ldots , a_m)$. The set of all weights is denoted by $\CW (a)=\CW (a; V)$. Let $\CA$ be the subalgebra of $\End_K(V)$ which is generated by the elements $a_1, \ldots , a_m$. Then $V$ is a (left) $\CA$-module. In this case, we also say that $V$ is an $a$-{\em module}. 
  The subspace of $V$,  $$V_\l :=V_\l (a):=\{ v\in V \, | \, a_iv=\l_iv\;\; {\rm  for }\;\; i=1, \ldots , m\},$$ is called the $\l$-{\em weight subspace} of $V$. Clearly, $V_\l$ is an $a$-module and the sum $\sum_{\l \in \CW (a)}V_\l = \bigoplus_{\l \in \CW (a)}V_\l$ is a direct sum of $a$-modules. If $$V= \bigoplus_{\l \in \CW (a)}V_\l$$ then we say that $V$ is a {\em weight} $a$-module. More generally, for each weight $\l\in \CW (a)$, the subspace $$V^\l :=V^\l (a):=\{ v\in V \, | \, (a_i-\l_i)^jv=0\;\; {\rm  for}\;\; j\gg 1\;\; {\rm  and }\;\; i=1, \ldots , m\},$$ is called the $\l$-{\em generalized weight subspace} of $V$. Clearly, $V_\l\subseteq V^\l$ and $V^\l$ is an $a$-module and the sum $\sum_{\l \in \CW (a)}V^\l = \bigoplus_{\l \in \CW (a)}V^\l$ is a direct sum of $a$-modules. If $$V= \bigoplus_{\l \in \CW (a)}V^\l$$ then we say that $V$ is a {\em generalized weight} $a$-module. If the field $K$ is an algebraically closed field and the vector space $V$ is a {\em locally finite} $a$-module, i.e. the linear maps $a_1, \ldots , a_m$ are locally finite (i.e. for all $v\in V$, $\dim_K(K[a_i]v)<\infty$ for $i=1, \ldots , m$) then the $a$-module $V$ is generalized weight module and vice versa. In particular, every finite dimensional vectors space $V$ is a generalized weight $a$-module.

For all elements $z\in Z(A_{x_1x_2})$, $i\in I_{A_{x_1x_2}}$ and $\beta \in \CR_{A_{x_1x_2}}$,
\begin{equation}\label{Zoox1x}
\o_{x_1}(zx_3^ix^\beta )=q^{-\beta_2}z(x_3-\alpha )^ix^\beta\;\;  \o_{x_2}(zx_3^ix^\beta )=q^{\beta_1}z(x_3-\mu )^ix^\beta\;\; {\rm and }\;\; \ad_{x_3}(zx_3^ix^\beta )=(\alpha \beta_1 +\mu \beta_2)zx_3^ix^\beta.
\end{equation}

\begin{proposition}\label{16Sep20}
Let $\o = (\o_{x_1}, \o_{x_2}, \ad_{x_3})$. Then the algebra $A_{x_1x_2}$ is a generalized weight $\o$-module,   $A_{x_1x_2}=\bigoplus_{\beta \in \CR}\bigg( \bigoplus_{i\in I} Z(A_{x_1x_2}) x_3^ix^\beta\bigg)$ is the direct sum of   generalized weight components where the sets $\CR=\CR_{A_{x_1x_2}}$ and $I=I_{A_{x_1x_2}}$ are  defined in Theorem \ref{ZD27Aug20},  $\CW (\o , A_{x_1x_2})=\{ (q^{-\beta_2}, q^{\beta_1}, \alpha \beta_1 +\mu \beta_2)\, | \, \beta =(\beta_1, \beta_2) \in \CR\}$, and  for each $\beta \in \CR$ the element $x^\beta$ is an $\o$-eigenvector with weight $(q^{-\beta_2}, q^{\beta_1}, \alpha \beta_1 +\mu \beta_2)$. 
\end{proposition}


{\it Proof}. By Theorem \ref{ZD27Aug20}, we have the direct sum as in the theorem. It remain to show that the vector spaces $V_\beta :=\bigoplus_{i\in I} Z(A_{x_1x_2}) x_3^ix^\beta$ are the generalized weight components. The theorem follows from the Claim.

{\sc Claim:} {\em The elements of the set 
$\CW':=\{ (q^{-\beta_2}, q^{\beta_1}, \alpha \beta_1 +\mu \beta_2)\, | \, \beta =(\beta_1, \beta_2) \in \CR\}$ are distinct.}

 In more detail, suppose that the Claim holds. Then, by (\ref{Zoox1x}),  $\CW (\o , A_{x_1x_2})\subseteq \CW'$. The reverse inclusion follows from the fact that for each $\beta \in \CR$, the element $x^\beta$ is a weight vector with weight $(q^{-\beta_2}, q^{\beta_1}, \alpha \beta_1 +\mu \beta_2)$, by (\ref{Zoox1x}). Hence, $\CW (\o , A_{x_1x_2})=\CW$. Now,  by the Claim and (\ref{Zoox1x}), the vector space $V_\beta$ is the generalized weight component that correspond to the weight $(q^{-\beta_2}, q^{\beta_1}, \alpha \beta_1 +\mu \beta_2)$. 
 
{\it Proof of the Claim.} Suppose that for  elements $\beta , \beta'\in \CR$, 
$$ (q^{-\beta_2}, q^{\beta_1}, \alpha \beta_1 +\mu \beta_2)=(q^{-\beta_2'}, q^{\beta_1'}, \alpha \beta_1' +\mu \beta_2').$$
we have to show that $\beta = \beta'$. This is obvious if $q$ is not root of unity, i.e. if the algebra $A_{x_1x_2}$ belongs to classes (C1)--(C3). 

So, we assume that $q$ is a primitive $n$'th root of unity, i.e.  the algebra $A_{x_1x_2}$ belongs to classes (C4)--(C11). Then the equality above is equivalent to the equalities:
\begin{eqnarray*}
  \beta_1'&\equiv & \beta_1 \mod n, \\ 
  \beta_2'&\equiv & \beta_2 \mod n,\\  
  (\beta_1'-\beta_1)\alpha &=& (\beta_2-\beta_2')\mu .
\end{eqnarray*}
{\em Cases (C4) and (C5):} In view of Lemma \ref{a16Sep20}.(1), it suffices to consider the case (C4). By Theorem \ref{ZD27Aug20}, $\CR = \Z\times \N_n$. Hence, $\beta_2'=\beta_2$. Then $\beta_1'=\beta_1$ since $\alpha =1$.

{\em Cases (C6) and (C7a):} In view of Lemma \ref{a16Sep20}.(2), it suffices to consider the case (C6).  The third equality, $\mu_2(\beta_1'-\beta_1)=-\mu_1(\beta_2-\beta_2')$, implies that $\beta_1'-\beta_1=k\mu_1$ and $\beta_2'-\beta_2=k\mu_2$ for some integer $k\in \Z$ (since $(\mu_1, \mu_2)=1$). Then $x^{\beta'}=x^\beta (x_1^{\mu_1}x_2^{\mu_2})^k$, and so 
 $\beta'=\beta$ since  $\CR = \mS' (\mu_1n, \mu_2n)$, by Theorem \ref{ZD27Aug20}. 
 
{\em Case (C7b):} Since $p=0$, $\alpha =1$ and $\mu\not\in \Q$, the third equality, $\beta_1'-\beta_1= (\beta_2-\beta_2')\mu $, yields the equality $\beta = \beta'$. 

{\em Cases (C8) and (C9):} In view of Lemma \ref{a16Sep20}.(1), it suffices to consider the case (C8). By Theorem \ref{ZD27Aug20}, $\CR = \N_{<pn}\times \N_{<n}$. The second equation, $\beta_2'\equiv \beta_2\mod j$, implies the equality $\beta_2'= \beta_2$. Then the third equality $\beta_1'-\beta_1=0$ (since $\alpha =1$) implies the equality $\beta=\beta'$. 
 
 {\em Case (C10):} Since $p>0$, $\alpha =1$ and $\mu\not\in \Fp$, the third equation, $\beta_1'-\beta_1 =(\beta_2-\beta_2')\mu$, yields the equivalences $$\beta_1'\equiv \beta_1 \mod p\;\; {\rm  and }\;\; \beta_2'\equiv \beta_2\mod p.$$  Then the first two equations yield $\beta_1'\equiv \beta_1 \mod pn$ and $\beta_2'\equiv \beta_2\mod pn$ since  $(p,n)=1$. Hence, $\beta =\beta'$ since   $\CR = \N_{<pn}\times \N_{<pn}$, by Theorem \ref{ZD27Aug20}.
 
 {\em Case (C11):} Recall that $p>0$, $\alpha =1$ and $\mu\in \Fp\backslash \{ 0, -1\}$. By Theorem \ref{ZD27Aug20}, $\CR = \N_{<pn}\times \N_{<n}$. Then the second equation, $\beta_2'\equiv \beta_2\mod n$, yields the equality $\beta_2'=\beta_2$. Then the third equation, $\beta_1'-\beta_1=0$, yields the equivalence  $\beta_1'\equiv \beta_1\mod p$. Then the first equation, $\beta_1'\equiv \beta_1\mod n$, implies that $\beta_1'\equiv \beta_1\mod pn$ since $(p,n)=1$. Hence, $\beta=\beta'$. $\Box $\\

{\bf The $\o$-eigenalgebra $\Es_{A_{x_1x_2}}(\o)$.} 
 By Theorem \ref{ZD27Aug20}, $$A_{x_1x_2}=\bigoplus_{i\in I_{A_{x_1x_2}}, \beta \in \CR_{A_{x_1x_2}}} Z(A_{x_1x_2}) x_3^ix^\beta.$$ So, each  element $a\in A_{x_1x_2}$ is a unique sum
 $a=\sum_{i\in I_{A_{x_1x_2}}, \beta \in \CR_{A_{x_1x_2}}} z_{i\beta} x_3^ix^\beta$ for some elements $z_{i\beta}\in Z(A_{x_1x_2})$. For $a\neq 0$, let 
 $${\rm Deg}_{x_3}(a):=\max \{ i \in I_{A_{x_1x_2}}\, | \, z_{i\beta}\neq 0\;\; {\rm for \; some} \;\ \beta \in \CR_{A_{x_1x_2}}\}. $$
Clearly, the sum 
$$\Es_{A_{x_1x_2}}(\o)=\bigoplus_{\l \in \CW (\o ,A_{x_1x_2})}A_{x_1x_2,_\l}$$ of all weight $\o$-submodules of $A_{x_1x_2}$ is the largest weight $\o$-submodule of the algebra $A_{x_1x_2}$. The space $\Es_{A_{x_1x_2}}(\o)$ is a subalgebra of $A_{x_1x_2}$ which is called the $\o$-{\em eigenalgebra} of $A_{x_1x_2}$. 

Corollary \ref{b16Sep20} is a description of the algebra $ \Es_{A_{x_1x_2}}(\o)$. It also shows that the algebra $A_{x_1x_2}$ is a free left/right $ \Es_{A_{x_1x_2}}(\o)$-module.  
 \begin{corollary}\label{b16Sep20}
\begin{enumerate}
\item $\Es_{A_{x_1x_2}}(\o)=\{a\in A_{x_1x_2} \, | \,{\rm Deg}_{x_3}(a)\leq 0\}= \bigoplus_{\beta \in \CR_{A_{x_1x_2}}}Z(A_{x_1x_2})x^\beta$; for all  $\l_\beta =(q^{-\beta_2}, q^{\beta_1}, \alpha \beta_1 +\mu \beta_2)\in \CW (\o , A_{x_1x_2})=\{ (q^{-\beta_2}, q^{\beta_1}, \alpha \beta_1 +\mu \beta_2)\, | \, \beta =(\beta_1, \beta_2) \in \CR_{A_{x_1x_2}}\}$, $A_{x_1x_2,\l_\beta}= Z(A_{x_1x_2})x^{\beta}$;  for all $z\in Z(A_{x_1x_2})$ and $\beta \in \CR_{A_{x_1x_2}}$, $\o_{x_1}(zx^\beta )=q^{-\beta_2}zx^\beta$, $\o_{x_2}(zx^\beta )=q^{\beta_1}zx^\beta$ and $\ad_{x_3}(zx^\beta )=(\alpha\beta_1+\mu\beta_2)zx^\beta$. 
\item $A_{x_1x_2}=\bigoplus_{i\in I_{A_{x_1x_2}}}x_3^i\Es_A(\o)=\bigoplus_{i\in I_{A_{x_1x_2}}}\Es_A(\o)x_3^i$ is a direct sum of $\Es_{A_{x_1x_2}}(\o)$-bimodules $x_3^i\Es_{A_{x_1x_2}}(\o)=\Es_{A_{x_1x_2}}(\o)x_3^i$. The algebra $A_{x_1x_2}$ is a free left/right $ \Es_{A_{x_1x_2}}(\o)$-module. 
\end{enumerate}
\end{corollary}


{\it Proof}. 1. Statement 1  follows from Proposition \ref{16Sep20} and (\ref{Zoox1x}).
In more detail, by (\ref{Zoox1x}), 
 $$\o_{x_1}(zx^\beta )=q^{-\beta_2}zx^\beta,\;\;  \o_{x_2}(zx^\beta )=q^{\beta_1}zx^\beta\;\; {\rm  and }\;\; \ad_{x_3}(zx^\beta )=(\alpha\beta_1+\mu\beta_2)zx^\beta.$$ If for some polynomial $\phi\in K[x_3]$,  $$ \o_{x_1}(\phi zx^\beta )=q^{-\beta_2}\phi  zx^\beta ,\;\; \o_{x_2}(\phi zx^\beta )=q^{\beta_1}\phi zx^\beta \;\; {\rm and} \; \; \ad_{x_3}(\phi zx^\beta )=(\alpha\beta_1+\mu\beta_2)\phi zx^\beta$$ then necessarily
 $\phi\in K[x_3]^\o\subseteq Z(A_{x_1x_2})$, and  statement 1 follows from Proposition \ref{16Sep20}.
 
 2. Statement 2 follows from statement 1 and Proposition \ref{16Sep20}.  $\Box $\\

 By Corollary \ref{b16Sep20}.(1), $\CW (\o , A_{x_1x_2})=\{ \l_\beta:=(q^{-\beta_2}, q^{\beta_1}, \alpha \beta_1 +\mu \beta_2)\, | \, \beta =(\beta_1, \beta_2) \in \CR_{A_{x_1x_2}}\}$ and  $A_{x_1x_2,\l_\beta}= Z(A_{x_1x_2})x^{\beta}$. For all $\beta , \beta'\in \CR_{A_{x_1x_2}}$, 
 $$ A_{x_1x_2,\l_\beta}A_{x_1x_2,\l_{\beta'}}= Z(A_{x_1x_2})x^{\beta +\beta'}=A_{x_1x_2,\l_{\beta +\beta'}}\;\; {\rm and}\;\; A_{x_1x_2,\l_{-\beta}}= Z(A_{x_1x_2})x^{-\beta}.$$
Therefore, the set $\CR_{x_1x_2}$ is an additive abelian group. By Corollary  \ref{a28Aug20} and Corollary \ref{a22Sep20}, the set $\CR_A$ is an additive submonoid of   the group $\CR_{x_1x_2}$. 
  Lemma \ref{a6Oct20} is an explicit description of the group $\CR_{x_1x_2}$.

 For each natural number $n$, let $Z_n=\Z/n\Z$. 
\begin{lemma}\label{a6Oct20}
 $$\CR_{A_{x_1x_2}}=\begin{cases}
\Z^2& \text{(C1): if  $q$ is not a root of unity and $p=0$},\\
\Z^2& \text{(C2): if  $q$ is not a root of unity, $p>0$ and} \\
 & (\alpha , \mu )\in \{ (1,\mu'), (0,1)\, | \, \mu'\in \mF_p\backslash \{ -1\}\} ,\\
\Z^2& \text{(C3): if  $q$ is not a root of unity, $p>0$, $\alpha =1$ and $\mu\not\in \mF_p $},\\
\Z\times\Z_{<n}& \text{(C4): if  $q$ is a primitive $n$'th  root of unity, $p=0$, $\alpha =1$ and $\mu =0$},\\
\Z_{<n}\times\Z& \text{(C5): if  $q$ is a primitive $n$'th  root of unity, $p=0$, $\alpha =0$ and $\mu =1$},\\
\Z^2/\Z (\mu_1n, \mu_2n)& \text{(C6): if  $q$ is a primitive $n$'th  root of unity, $p=0$, $\alpha =1$ and $\mu =-\frac{\mu_1}{\mu_2}$},\\
& \text{ where $\mu_1,\mu_2\in \N_+$ such that $\mu_1\neq \mu_2$ and $(\mu_1,\mu_2)=1$},\\
\Z^2/\Z(-\mu_1n, \mu_2n) &\text{(C7a): if  $q$ is a primitive $n$'th  root of unity, $p=0$, $\alpha =1$ and $\mu =\frac{\mu_1}{\mu_2}$},\\
& \text{ where $\mu_1,\mu_2\in \N_+$ such that $(\mu_1,\mu_2)=1$},\\
\Z^2& \text{(C7b): if  $q$ is a primitive $n$'th  root of unity, $p=0$, $\alpha =1$ and $\mu\not\in \Q$},\\
\Z_{<pn}\times\Z_{<n}& \text{(C8): if  $q$ is  a primitive $n$'th  root of unity, $p>0$, $\alpha =1$, $\mu =0$ and $p\nmid n$}, \\
\Z_{<n}\times\Z_{<pn}& \text{C(9): if  $q$ is  a primitive $n$'th  root of unity, $p>0$, $\alpha =0$, $\mu =1$ and $p\nmid n$,} \\
\Z_{<pn}\times\Z_{<pn}& \text{(C10): if  $q$ is  a primitive $n$'th  root of unity, $p>0$, $\alpha =1$, $\mu\not\in \mF_p $  and $p\nmid n$},\\
\Z_{<pn}\times\Z_{<n}& \text{(C11): if  $q$ is  a primitive $n$'th  root of unity, $p>0$, $\alpha =1$}, \\
& \text{$\mu\in \mF_p \backslash \{ 0,-1\}=\{ 1,\ldots , p-2\}$  and $p\nmid n$.}
\end{cases}
$$
\end{lemma} 
 
  {\it Proof}. The lemma follows from Theorem \ref{ZD27Aug20}. $\Box$ 
 
 By Corollary \ref{b16Sep20}.(1), $\CW (\o , A_{x_1x_2})=\{ \l_\beta=(q^{-\beta_2}, q^{\beta_1}, \alpha \beta_1 +\mu \beta_2)\, | \, \beta =(\beta_1, \beta_2) \in \CR_{A_{x_1x_2}}\}$, i.e. the map $\CR_{A_{x_1x_2}}\ra \CW (\o , A_{x_1x_2})$, $\beta \mapsto \l_\beta $ is a bijection. By Lemma \ref{a6Oct20}, the set  $\CW (\o , A_{x_1x_2})$ is an abelian group isomorphic to $\CR_{A_{x_1x_2}}$ via the bijection above. So, the group structure on the set $\CW (\o , A_{x_1x_2})$ is given by the rule:  for all elements $\beta , \beta'\in \CR_{A_{x_1x_2}}$, $\l_\beta +\l_{\beta'}:= \l_{\beta +\beta'}$. \\
 
{\bf Simplicity criterion for the algebra $A_{x_1x_2}$.} 
\begin{theorem}\label{A14Sep20}
The algebra $A_{x_1x_2}$ is a simple algebra iff the algebra $A$ belongs to the case (C1) and  or (C7b) of Theorem \ref{C24Sep19} iff $Z(A_{x_1x_2)}=K$. 
\end{theorem}

{\it Proof}. The second `iff' is obvious, see Theorem \ref{AC24Sep19}. 

Let us prove that the first `iff' holds.

$(\Rightarrow )$ The centre of a simple algebra is a field, and so the implication follows from Theorem \ref{AC24Sep19}. 

$(\Leftarrow )$ Let $J$ be a nonzero ideal of the algebra $A_{x_1x_2}$. We have to show that $J=(1)$.

Suppose that the algebra $A$ belongs to the case (C7b). Then the elements  $\{x^\beta\}_{\beta \in \Z^2}$ are eigenvectors of the derivation $\ad_{x_3}$ of the algebra $A_{x_1x_2}$ with distinct eigenvalues (since $\alpha = 1\in \Q$ and $\mu \not\in \Q$). Using the $\Z^2$-grading of the algebra $A_{x_1x_2}$, we conclude that $\v x^\beta\in J$ for some nonzero element $\v\in K[x_3]$ and $\beta \in \Z^2$. Hence, $\v \in J$ since the element $x^\beta$ is a unit.
 We may assume that the element $\v$ has the least possible degree $d$ in the variable $x_3$. If $d=0$ then $\v\in K^\times$ and we are done. Suppose that $d>0$, we seek a contradiction. Recall that $K[x_3]^G=K$ (Proposition \ref{E24Sep19}). The element $$0\neq (1-\o_{x_1})(\v (x_3))=v(x_3)-v(x_3-1 )\in J $$  has degree $d-1$, a contradiction. 

Suppose that the algebra $A$ belongs to case (C1). Then the elements  $\{x^\beta\}_{\beta \in \Z^2}$ are common eigenvectors of the commuting inner automorphisms $\o_{x_1}$ and $\o_{x_2}$ of the algebra $A_{x_1x_2}$ with distinct weights (since $q$ is not a root of unity). Furthermore, the algebra $A=\bigoplus_{\beta \in \Z^2} K[x_3]x^\beta$ is a direct sum of generalized weight subspaces $K[x_3]x^\beta$. Hence,  $\v x^\beta\in J$ for some nonzero element $\v\in K[x_3]$ and $\beta \in \Z^2$. Hence, $\v \in J$ since the element $x^\beta$ is a unit.
 We may assume that the element $\v$ has the least possible degree $d$ in the variable $x_3$. If $d=0$ then $\v\in K^\times$ and we are done. Suppose that $d>0$, we seek a contradiction. Recall that $K[x_3]^G=K$ (Proposition \ref{E24Sep19}). Hence, one of the elements  $$(1-\o_{x_1})(\v (x_3))=v(x_3)-v(x_3-\alpha )\in J\;\; {\rm  or} \;\; (1-\o_{x_2})(\v (x_3))=v(x_3)-v(x_3-\mu )\in J$$ is nonzero and has degree strictly smaller than $d$, a contradiction. 
 $\Box $\\

{\bf The ideal structure and the  prime spectrum of the algebra $A_{x_1x_2}$.} For an algebra $R$, we denote by $\CI(R)$ the set of ideals of $R$. 

\begin{theorem}\label{15Sep20}
The map $$\CI (Z(A_{x_1x_2}))\ra \CI (A_{x_1x_2}), \;\; \gp \mapsto A_{x_1x_2}\gp$$ is a bijection with the inverse $Q\mapsto Z(A_{x_1x_2})\cap Q$.
\end{theorem}

{\it Proof}. In general,  the theorem is true for all simple algebras. In particular, it is true  for simple algebras $A_{x_1x_2}$, i.e. the algebras  that belong to the cases (C1) and (C7b) (Theorem \ref{A14Sep20}).

So, we may assume that the algebra $A_{x_1x_2}$ does not belong to  the cases (C1) and (C7b). In particular,  the centre of the algebra $A_{x_1x_2}$ is not the field $K$ (Theorem \ref{AC24Sep19}). By Lemma \ref{a16Sep20} and the comments before it, it suffices to consider the cases (C2)--(C4), (C6), (C8), (C10) and (C11). So, we may assume that the algebra $A_{x_1x_2}$ belongs to one of these cases. Let $J $ be a nonzero ideal of the algebra $A_{x_1x_2}$. By Theorem \ref{ZD27Aug20}, in order to finish the prove of the theorem we have to show that 
$$J=J\cap Z(A_{x_1x_2}).$$
By Theorem \ref{ZD27Aug20}, each element $a\in J$ is a unique sum $a=\sum_{i\in I , \beta \in \CR}z_{i\beta}x_3^ix^\beta$ where $z_{i\beta }\in Z(A_{x_1x_2})$. So, we have to show that all elements $z_{i\beta}$ belong to the ideal $J$. The ideal $J$ is an  $\o$-module. By Proposition \ref{16Sep20}, all the elements $\sum_{i\in I}z_{i\beta}x_3^ix^\beta$, where $ \beta \in \CR$, belong to $J$. So, we may assume that $a=\sum_{i\in I}z_{i\beta}x_3^ix^\beta$ for some $ \beta \in \CR$. The element $x^\beta$ is a unit element of the algebra $A_{x_1x_2}$. Hence, we may assume that $\beta =0$ (by dividing the element $a$ on the right by $x^\beta$), i.e. 
$$ a=a(x_3)=\sum_{i=0}^d z_{i0}x_3^i$$
for some $d=\Deg_{x_3}\geq 0$ such that $z_{d0}\neq 0$. The case when $D=0$ is obvious. So, we assume that $d>0$ and will use induction on $d$.

(i) {\em Suppose that the algebra $A_{x_1x_2}$ belongs to one of the cases (C4) or (C6) (resp.,(C2), (C8) or (C11)). } In each of these cases, $p=0$, $I=\N$ and $\alpha =1$ (resp., $p>0$, $I=\N_{<p}$, $\alpha =1$ and $d<p$). The element
$$J\ni (1-\o_{x_1})(a)=a(x_3)-a(x_3-1)=dz_{d0}x_3^{d-1}+\cdots $$
has degree $d-1$. By induction $dz_{d0}\in \ga$, and so $z_{d0}\in \ga$ since the element $d$ is a unit in $A_{x_1x_2}$ in all five cases. The element $a-z_{d0}x^d\in J$ has degree $<d$. By induction, all coefficients $z_{i0}$ belong to $J$. 

(ii) {\em Suppose that the algebra $A_{x_1x_2}$ belongs to one of the cases (C3),  or (C10)}. In each of these cases, $p>0$, $I=\N_{<p^2}$,  $\alpha =1$ and $\mu \not\in \Fp$. We may assume that $d\geq p$ otherwise we repeat word for word the proof of the statement (i). In this case the element $a $ is a unique sum $\sum_{j=0}^{p-1}\psi_jx_3^j$ for some elements $\psi_j\in \sum_{i=0}^{p-1}\xi_{ji}(x_3^p-x_3)^i$ where $\xi_{ji}\in K[x_3]^G\subseteq Z(A_{x_1x_2})$. We have to show that all the elements $\xi_{ji}$ belong to the ideal $J$.  Let 
$$D'(a ):=\max\{j \, | \, \psi_j\neq 0\}.$$
Suppose that $D'(a )=0$, i.e. $a =\psi_0=\sum_{i=0}^{p-1}\xi_{0i}(x_3^p-x_3)^i=\sum_{i=0}^{\d}\xi_{0i}(x_3^p-x_3)^i$ for some natural number $\d<p$ such that $\xi_{0\d}\neq 0$. Clearly, $\d \geq 1$ as $d\geq p$. Then $\mu^p-\mu \neq 0$ (since $\mu \not\in \Fp )$ and 
\begin{eqnarray*}
 J&\ni & (1-\o_{x_2})(a (x_3))=a(x_3)-a(x_3-\mu )=\sum_{i=0}^{\d}\xi_{0i}(x_3^p-x_3)^i-
\sum_{i=0}^{\d}\xi_{0i}(x_3^p-x_3-(\mu^p-\mu))^i \\
 &=& \d (\mu^p-\mu )\xi_{0\d}(x_3^p-x_3)^{\d -1}+\cdots .
\end{eqnarray*}
Then by induction on $\d$,  $\xi_{0i}\in J$ for all $i$ such that  $1\leq i\leq \d$. Now, $\xi_{00}=v-\sum_{i=1}^{\d} \xi_{0i}(x_3^p-x_3)^i\in J$. So, all the  elements $\xi_{0i}$ belong to $J$.

 Suppose that $D'=D'(a )>0$. Now,  
$$J\ni (1-\o_{x_1})(a (x_3))=a(x_3)-a(x_3-1 )=\sum_{j=0}^{D'}\psi_j(x_3^j-(x_3-1)^j)=D'\psi_{D'}x_3^{D'-1}+\cdots .
$$
By induction on $D'$,  $\xi_{D'i}\in J$. The element $a':=a-\psi_{D'}x_3^{D'}\in J$ has $D'(a')<D'(a)=D'$. By induction on $D'$, all the coefficients $\xi_{ji}\in J$, as required. This finishes the proof of  statement (ii) and the theorem. $\Box$

Theorem \ref{AA15Sep20} is a classification of prime ideals of the algebra $A_{x_1x_2}$. 

\begin{theorem}\label{AA15Sep20}
The map $$\Spec (Z(A_{x_1x_2}))\ra \Spec (A_{x_1x_2}), \;\; \gp \mapsto A_{x_1x_2}\gp$$ is a bijection with the inverse $Q\mapsto Z(A_{x_1x_2})\cap Q$. In particular, every nonzero prime ideal of $A_{x_1x_2}$ meets the centre.
\end{theorem}

{\it Proof.} The theorem follows at once from Theorem \ref{15Sep20}. $\Box$


\section{The prime spectrum of the algebra $A$}\label{SPECA}

The aim of the section is to classify the  prime ideals of the algebra  $A$ (see (\ref{SpAAa}), (\ref{SpAAa3}) and   Theorem \ref{18Sep20}).\\

For an algebra $R$, let $\Spec\,(R)$ be the set of its prime ideals. The set $(\Spec\,(R), \subseteq)$ is a partially ordered set (poset) with respect to inclusion of prime ideals. A prime ideal $\gp$ of an algebra $R$ is called a \emph{completely prime ideal} if $R/\gp$ is a domain. We denote by  $\Spec_c(R)$ the set of completely prime ideals of $R$; it is called the \emph{completely prime spectrum} of $R$. The annihilator of a simple $R$-module is a prime ideal of $R$. That kind of prime ideals are called {\em primitive}. The set of all primitive ideal of $R$ is denoted by ${\rm Prim} (R)$. By the very definition, the poset $(\Spec\,(R), \subseteq)$ is determined by the containment information about prime ideals of $R$. In general, it is difficult to say whether one prime ideal contains or does not contain another prime ideal. Proposition \ref{aA12Mar15}.(3) is a general result about containments of primes. 
It splits the prime spectrum of a ring into two disjoint subsets such that elements of one are not contained in elements of the other. 

 Each element $r \in R$ determines two maps from $R$ to $R$,  $r \cdot: x \mapsto rx$ and $\cdot r: x \mapsto xr$ where $x \in R.$ An element $a \in R$ is called a \emph{normal element} if $Ra=aR$.
\begin{proposition} \label{aA12Mar15} 
(\cite{Bav-Lu-BL-qAge}.)  
Let $R$ be a Noetherian ring and $s$ be an element of $R$ such that $S_s := \{s^i \,|\, i\in \N \}$ is a left denominator set  of the ring $R$ and $(s^i)= (s)^i$ for all $i\geqslant 1$ (e.g., $s$ is a normal element such that $\ker(\cdot s_{R}) \subseteq \ker(s_R \cdot)$). Then
$\Spec\,(R)= \Spec(R, s) \, \sqcup \, \Spec_s(R) $ where $\Spec(R,s):= \{ \gp \in \Spec\,(R)\,|\, s \in \gp \}$,  $\Spec_s(R)= \{\gq \in \Spec\,(R)\, |\, s \notin \gq \}$ and
\begin{enumerate}
\item  the map $\Spec\,(R, s) \rightarrow \Spec\,(R/(s)),\,\, \gp \mapsto \gp/(s)$, is a bijection with the inverse $\gq \mapsto \pi^{-1}(\gq)$ where $\pi: R \rightarrow R/(s), r \mapsto r+(s),$
\item the map $\Spec_s(R) \rightarrow \Spec\,(R_s), \,\, \gp \mapsto S_s^{-1}\gp,$ is a bijection with the inverse $\gq \mapsto \sigma^{-1}(\gq)$ where $\sigma: R \rightarrow R_s:= S_s^{-1}R, \,r \mapsto \frac{r}{1}$. 
\item For all $\gp \in \Spec\,(R, s)$ and $\gq \in \Spec_s(R), \gp \not\subseteq \gq.$
\end{enumerate}
\end{proposition}
{\bf The prime spectra and simple modules for the algebras $A/(x_1)$ and $A/(x_2)$.} 
\begin{eqnarray*}
A/(x_1)&=&K\langle x_2,x_3\, | \, x_2x_3=(x_3-\mu )x_2\rangle \simeq \begin{cases}
\mA & \text{if }\mu\neq 0,\\
K[x_2,x_3]& \text{if }\mu =0.\\
\end{cases}\\
A/(x_2)&=&K\langle x_1,x_3\, | \, x_1x_3=(x_3-\alpha )x_1\rangle \simeq \begin{cases}
\mA & \text{if }\alpha\neq 0,\\
K[x_1,x_3]& \text{if }\alpha =0.\\
\end{cases}
\end{eqnarray*}
The algebras $A/(x_1)$ and $A/(x_2)$ are {\em subalgebras}  of $A$. The prime ideals   and simple modules for the algebra $\mA$   are classified in \cite{SpecWeylcharp}.

For each subset $\ga$ of $A$, $V(\ga ) :=\{ \gp \in \Spec (A) \, | \, \ga \subseteq \gp \}$. The element $x_1x_2\in A$ is a regular normal element. Therefore, $A\subseteq A_{x_1x_2}$. By Proposition \ref{aA12Mar15}, 
\begin{equation}\label{SpAAa}
\Spec (A)= \bigg( V(x_1)\cup V(x_2)\bigg)\coprod \Spec_{x_1x_2}(A)
\end{equation}
since $\Spec (A, x_1x_2)= V(x_1x_2)=V(x_1)\cup V(x_2)$
 as the elements $x_1$ and $x_2$ are normal,  and the map 
\begin{equation}\label{SpAAa1}
\Spec_{x_1x_2}(A)\ra \Spec (A_{x_1x_2}), \;\; P\mapsto P_{x_1x_2}
\end{equation}
is a {\em bijection with the inverse} $Q\mapsto A\cap Q$. By Theorem \ref{AA15Sep20}, the map  $$\Spec (Z(A_{x_1x_2}))\ra \Spec (A_{x_1x_2}), \;\; \gp \mapsto A_{x_1x_2}\gp$$ is a bijection with the inverse $Q\mapsto Z(A_{x_1x_2})\cap Q$. Therefore, the map
\begin{equation}\label{SpAAa2}
\Spec (A_{x_1x_2})\ra \Spec_{x_1x_2}(A), \;\; \gq \mapsto A\cap \gq A_{x_1x_2}
\end{equation}
is a bijection with the inverse $P\mapsto Z(A_{x_1x_2})\cap P_{x_1x_2}$. Since the elements $x_1$ and $x_2$ are regular and normal elements of the algebra $A$ we have the equalities
\begin{equation}\label{SpAAa3}
V(x_1)=\Spec (A/(x_2))\;\; {\rm and}\;\; V(x_2)=\Spec (A/(x_1)).
\end{equation}

{\bf  The prime spectrum of  $A$: Cases (C1) and (C7b).}

\begin{theorem}\label{SpecC1C7b}
The following statements are equivalent:
\begin{enumerate}
\item The algebra  $A$ belongs to one of the cases (C1) or (C7b). 
\item $\Spec (A)=\{ 0\} \cup V(x_1)\cup V(x_2)=\{ 0\} \cup \Spec (A/(x_1))\cup \Spec (A/(x_2))$.
\item $\Spec_{x_1x_2} (A)=\{ 0\}$.
\item $\Spec (A_{x_1x_2})=\{ 0\}$.
\item The algebra $A_{x_1x_2}$ is a simple algebra.
\end{enumerate}
\end{theorem}

{\it Proof}. By (\ref{SpAAa}), $\Spec (A)=\{ 0\} \cup V(x_1)\cup V(x_2)$ iff  $\Spec_{x_1x_2} (A)=\{ 0\}$ iff  $\Spec (A_{x_1x_2})=\{ 0\}$ iff the algebra $A_{x_1x_2}$ is a simple algebra iff the algebra  $A$ belongs to one of the cases (C1) or (C7b), by Theorem \ref{A14Sep20}.  $\Box$ 

The set $T:=T_A:=Z(A)\cap \{ x^\beta\, | \, \beta \in \N^2\}$ is an Ore set of the algebra $A$ that belongs to the centre of $A$. Lemma \ref{a15Sep20} describes generators of the monoid $T$.
\begin{lemma}\label{a15Sep20}
Let $M_{\S} :=\langle x_1^{pn}, x_1^{\xi_i n}x_2^{in},  x_2^{pn}\rangle_{i=1, \ldots ,p-1}$ be the multiplicative  submonoid of $\S$ (which is generated by the elements in the brackets)  where $\xi_i$ is a unique natural number such that $0\leq \xi_i \leq p-1$ and $ \xi_i\equiv -i\mu \mod p$ (see Theorem \ref{C24Sep19}) and $T=Z(A)\cap \{ x^\beta\, | \, \beta \in \N^2\}$. Then
$$T=\begin{cases}
\{ 1\} & \text{(C1): if  $q$ is not a root of unity and $p=0$},\\
\{ 1\} & \text{(C2): if  $q$ is not a root of unity, $p>0$ and} \\
 & (\alpha , \mu )\in \{ (1,\mu'), (0,1)\, | \, \mu'\in \mF_p\backslash \{ -1\}\} ,\\
\{ 1\} & \text{(C3): if  $q$ is not a root of unity, $p>0$, $\alpha =1$ and $\mu\not\in \mF_p $},\\
\langle x_2^n\rangle & \text{(C4): if  $q$ is a primitive $n$'th  root of unity, $p=0$, $\alpha =1$ and $\mu =0$},\\
\langle x_1^n\rangle & \text{(C5): if  $q$ is a primitive $n$'th  root of unity, $p=0$, $\alpha =0$ and $\mu =1$},\\
\langle x_1^{\mu_1n}x_2\rangle & \text{(C6): if  $q$ is a primitive $n$'th  root of unity, $p=0$, $\alpha =1$ and $\mu =-\frac{\mu_1}{\mu_2}$},\\
& \text{ where $\mu_1,\mu_2\in \N_+$ such that $\mu_1\neq \mu_2$ and $(\mu_1,\mu_2)=1$},\\
\{ 1\} & \text{(C7a): if  $q$ is a primitive $n$'th  root of unity, $p=0$, $\alpha =1$ and $\mu =\frac{\mu_1}{\mu_2}$},\\
& \text{ where $\mu_1,\mu_2\in \N_+$ such that  $(\mu_1,\mu_2)=1$},\\
\{ 1\} & \text{(C7b): if  $q$ is a primitive $n$'th  root of unity, $p=0$, $\alpha =1$ and $\mu\not\in \Q $},\\
\langle x_1^{pn}, x_2^{n}\rangle & \text{(C8): if  $q$ is  a primitive $n$'th  root of unity, $p>0$, $\alpha =1$, $\mu =0$ and $p\nmid n$}, \\
\langle x_1^{ n}, x_2^{ pn}\rangle & \text{(C9): if  $q$ is  a primitive $n$'th  root of unity, $p>0$, $\alpha =0$, $\mu =1$ and $p\nmid n$,} \\
\langle x_1^{pn}, x_2^{ pn}\rangle& \text{(C10): if  $q$ is  a primitive $n$'th  root of unity, $p>0$, $\alpha =1$, $\mu\not\in \mF_p $  and $p\nmid n$},\\
M_{\S} & \text{(C11): if  $q$ is  a primitive $n$'th  root of unity, $p>0$, $\alpha =1$}, \\
& \text{$\mu\in \mF_p \backslash \{ 0,-1\}=\{ 1,\ldots , p-2\}$  and $p\nmid n$.} \\
\end{cases}
$$
\end{lemma}

{\it Proof}. The lemma follows from the description of the centre of the algebra $A_{x_1x_2}$  (Theorem \ref{AC24Sep19}. $\Box $

\begin{proposition}\label{A15Sep20}
Let $T=Z(A)\cap \{ x^\beta\, | \, \beta \in \N^2\}$. Then the following statements are equivalent:
\begin{enumerate}
\item The algebra  $A$ is
not of type (C7a). 
\item $Z(A_{x_1x_2})=T^{-1}Z(A)$.
\end{enumerate}
\end{proposition}
{\it Proof}. Suppose that the algebra $A$ is of type   (C7a). Then $T=\{ 1\}$, and so 
 $T^{-1}Z(A)=Z(A)=K\neq Z(A_{x_1x_2})$, by Theorem \ref{C24Sep19} and Theorem \ref{AC24Sep19}.

Suppose that the algebra $A$ is not  of type (C7a). We have to show that 
 $Z(A_{x_1x_2})=T^{-1}Z(A)$. This follows straightaway from the description of the set $T$ (Lemma \ref{a15Sep20}) and the descriptions of the centres of the algebras $A$ (Theorem \ref{C24Sep19})  and $A_{x_1x_2}$ (Theorem \ref{AC24Sep19}). $\Box$

\begin{lemma}\label{XA15Sep20}
Let $T=Z(A)\cap \{ x^\beta\, | \, \beta \in \N^2\}$. Then the following statements are equivalent:
\begin{enumerate}
\item The algebra  $A$ is of type (C6) or (C8)--(C11).
\item $A_{x_1x_2}=T^{-1}A$.  
\item $Z(A_{x_1x_2})=T^{-1}Z(A)$.
\end{enumerate}
\end{lemma}

{\em Remark.}  Lemma \ref{XA15Sep20} states that if the algebra $A$ is  of type (C6) or (C8)--(C11) then the localization $A_{x_1x_2}$ of the algebra $A$ at  the  Ore set that is generated by a {\em normal} element $x_1x_2$ is in fact the localization of the algebra $A$ at an Ore set that is generated by {\em central} elements.  Also the centre of the algebra $A_{x_1x_2}$ is the localization of the centre of the algebra $A$ at  an Ore set that is generated by {\em central} elements. In general, both these statements are badly wrong (as the algebras $A$ of type (C7a) show:  $T_A=\{1\}$, $Z(A_{x_1x_2})=K[(x_1^{\mu_1n}x_2^{-\mu_2n})^{\pm 1} ]\neq K=Z(A)=T^{-1}_AZ(A)$ and  $A_{x_1x_2}\neq A=T^{-1}_AA$).

{\it Proof.} The lemma follows from Lemma \ref{a15Sep20} and Proposition \ref{A15Sep20}. $\Box$

Theorem \ref{18Sep20} is an explicit description of the set $\Spec_{x_1x_2}(A)$. 

\begin{theorem}\label{18Sep20}

\begin{enumerate}
\item Suppose that the algebra $A$ is not of type (C7a) or (C11). Then the map 
$$ \Spec (Z(A))\backslash V(T_A)\ra \Spec_{x_1x_2}(A), \;\; \gp \mapsto \gp A=\bigoplus_{i\in I_A, \beta \in \CR_A}\gp x_3^ix^\beta$$ is a bijection with the inverse $P\mapsto P\cap Z(A)$  where the set $T_A$ is  defined in Lemma  \ref{a15Sep20} and  the sets $I_A$ and $\CR_A$ are defined in Theorem \ref{D27Aug20}. 
\item Suppose that the algebra $A$ is of type (C7a) or (C11). Then the map 
$$ \Spec (Z(A_{x_1x_2}))\ra \Spec_{x_1x_2}(A), \;\; \gq \mapsto A\cap \gq A_{x_1x_2}$$ is a bijection with the inverse $P\mapsto P_{x_1x_2}\cap Z(A_{x_1x_2})$. 
\end{enumerate}
\end{theorem}
 
{\it Proof}. 1. Let $T=T_A$. By Proposition \ref{A15Sep20},  $Z(A_{x_1x_2})=T^{-1}Z(A)$, and so $\Spec (Z(A_{x_1x_2}))=\Spec (Z(A))\backslash V(T)$. Now, by (\ref{SpAAa2}), the map 
$$\Spec (Z(A))\backslash V(T)\ra \Spec_{x_1x_2}(A), \;\; \gp \mapsto A\cap \gp A_{x_1x_2}$$ is a bijection with the inverse $P\mapsto P_{x_1x_2}\cap Z(A_{x_1x_2})$. By Proposition \ref{A15Sep20} and Theorem  \ref{D27Aug20},  
$$A\cap \gp A_{x_1x_2}=A\cap \gp T^{-1}A=\bigoplus_{i\in I_A, \beta \in \CR_A}\Big(Z(A) \cap \gp T^{-1}Z(A)\Big) x_3^ix^\beta=\bigoplus_{i\in I_A, \beta \in \CR_A}\gp x_3^ix^\beta,$$
 $\gp A_{x_1x_2}\cap Z(A_{x_1x_2})=\gp A_{x_1x_2}\cap T^{-1}Z(A)=\gp T^{-1}Z(A)=T^{-1}\gp$, and statement 1 follows. 

2. Statement 2 is (\ref{SpAAa2}). $\Box $

In contrast to Theorem \ref{AA15Sep20} not every nonzero prime ideal of the set $\Spec_{x_1x_2}(A)$ meets the centre of the algebra $A$. 

\begin{corollary}\label{a18Sep20}
The following statements are equivalent:
\begin{enumerate}
\item For all nonzero prime ideals $P\in\Spec_{x_1x_2}(A)$, $P\cap    Z(A)\neq 0$.
\item The algebra $A$ is not of the type (C7a).
\end{enumerate}
\end{corollary}

{\it Proof}. By Proposition \ref{A15Sep20},  
 if the algebra $A$ is not of the type (C7a) then 
$P\cap    Z(A)\neq 0$ for all nonzero prime ideals $P\in\Spec_{A_{x_1x_2}}(A)$ (since the algebra $A$ is an essential left/right $A$-submodule of $A_{x_1x_2}$, $Z(A_{x_1x_2})= T_A^{-1}Z(A)$ and Theorem \ref{AA15Sep20}).

If the algebra $A$ is of the type (C7a) then $Z(A_{x_1x_2})\neq K=Z(A)$, and so $P\cap    Z(A)=P\cap K= 0$ for all nonzero prime ideals $P\in\Spec_{x_1x_2}(A)$.  $\Box $

\begin{corollary}\label{b18Sep20}
Suppose that the algebra $A$ belongs to the case (C7a). Then  for all nonzero prime ideals $P\in\Spec_{x_1x_2}(A)$, $P\cap    Z(A)=0$.
\end{corollary}

{\it Proof}. See the proof of Corollary \ref{a18Sep20}. $\Box$

For a ring $R$ and its prime ideal $P$, $\cht (P):=\Kdim (R/P)$ is called the {\em co-height} of $P$.

The algebra $A$ is an $\o$-submodule of the algebra $A_{x_1x_2}$.

\begin{lemma}\label{b14Sep20}
Every prime ideal  $P\in \Spec_{x_1x_2}(A)$ is an $\o$-submodule of the algebra $A$. 
\end{lemma}

{\it Proof}.  Since  $P\in \Spec_{x_1x_2}(A)$, $x_i\not\in P$ for $i=1,2$. Then it follows from the inclusions:
\begin{eqnarray*}
P\supseteq (x_i)P&=&(x_iP)=(\o_{x_i}(P)x_i)
=\o_{x_i}(P)(x_i), \\
P\supseteq P(x_i)&=&(Px_i)=(x_i\o_{x_i}^{-1}(P))
=(x_i)\o_{x_i}^{-1}(P),
\end{eqnarray*}
that $ \o_{x_i}^{\pm 1}(P)\subseteq P$, and so $\o_{x_i}(P)=P$. Clealry, $\ad_{x_3}(P)\subseteq P$, and the lemma follows. $\Box $

\begin{corollary}\label{c18Sep20}
For all nonzero prime ideals  $P\in \Spec (A)$, $P\cap \Es_A (\o )\neq 0$. 
\end{corollary} 

{\it Proof}. The corollary holds for all nonzero  prime ideals of the algebra $A$ that contains the element $x_1x_2\in \Es_A(\o )$.
 By (\ref{SpAAa}), we may assume that $0\neq P\in \Spec_{x_1x_2}(A)$.

If the algebra $A$ is not the type (C7a) then the corollary follows from Corollary \ref{a18Sep20}.

Suppose that the algebra is  of the type (C7a). 
Since $Z(A_{x_1x_2})=K[x_1^{\mu_1n}x_2^{-\mu_2n})^{\pm 1}]$ (Theorem \ref{AC24Sep19}), the corollary follows from Theorem \ref{18Sep20}.(2). $\Box$


\section{The automorphism group of the algebra $A$}\label{AUTTH5.1.(1)}

 The aim of the section is to describe the automorphism group $\Aut_K(A)$ (Theorem \ref{A27Aug20}),  to classify $\o$-weight vectors  (Theorem \ref{20Sep20}) and normal elements  (Proposition  \ref{A20Sep20}) of the algebra $A$. Furthermore, we will show that every normal element of the algebra $A$ is an $\o$-weight vector, and vice versa (Proposition \ref{A20Sep20}). \\

{\bf The automorphism group $\Aut_K(A)$.}  Proposition \ref{a21Sep20} is the key fact in finding an explicit description  for the automorphism group of the algebra $A$ (Theorem \ref{A27Aug20}).

\begin{proposition}\label{a21Sep20}

\begin{enumerate}
\item The ideals $(x_1)$ and $(x_2)$ of the algebra $A$ has co-height 2. The ideals $(x_1)$ and $(x_2)$ of the algebra $A$ are the only ideals that have co-height 2 iff the algebra $A$ 
is of type (C1)--(C7).
\item If the algebra $A$ is of type (C8) (resp., (C9)) then the ideal $(x_1)$ (resp., $(x_2))$ is the only height 1 prime ideal $P$ of $A$  such that the algebra $A/P$ is commutative.

\item Suppose that the algebra $A$ is of type (C10) or (C11). Then  the ideals $(x_1)$ and $(x_2)$ of the algebra $A$ is the only pair of prime ideals  of the algebra $A$ of  height 1, $P$ and $Q$,    the factor algebra $A/(P+Q)$ is commutative domain. 
\end{enumerate}
\end{proposition}

{\it Proof}. 1. Clearly, $\cht (x_i)=\Kdim (A/(x_i))=2$. The second part of the corrollary follows from the classification of prime  ideals in the set $\Spec_{x_1x_2}(A)$ (Theorem  \ref{18Sep20}) and the descriptions of the centres of the algebras $A$ and $A_{x_1x_2}$ (Theorem \ref{C24Sep19} and Theorem \ref{AC24Sep19}). 

2.  In the view of the $(x_1,x_2)$-symmetry, it suffices to consider the case (C8). The algebra $A/(x_1)=K[x_2,x_3]$ is obviously a commutative domain and the algebra $A/(x_2)$ 
is obviously not (since $x_1x_3=(x_3-1)x_1$). 
 
 Suppose that $P\in \Spec_{x_1x_2}(A)$. By Theorem \ref{D27Aug20}, $A=\bigoplus_{i\in I, \beta \in \CR}Z(A)x_3^ix^\beta$.    Then the algebra $$A/(P)=\bigoplus_{i\in I, \beta \in \CR}\Big(Z(A)/(P)\Big)x_3^ix^\beta$$ is not commutative, by Theorem \ref{D27Aug20}, and the result follows (from the classification of prime ideals of the algebra $A$). 
 
3. The algebra $A/(x_1,x_2)=K[x_3]$ is obviously a commutative domain. 

Suppose that $P, Q\in \Spec_{x_1x_2}(A)$. By Theorem \ref{ZD27Aug20}, $A_{x_1x_2}=\bigoplus_{i\in I, \beta \in \CR}Z(A_{x_1x_2})x_3^ix^\beta$.    Then the algebra $$\Big(A/(P,Q)\Big)_{x_1x_2}=A_{x_1x_2}/(P,Q)=\bigoplus_{i\in I, \beta \in \CR}\Big(Z(A_{x_1x_2})/(P,Q)\Big)x_3^ix^\beta$$ is not commutative, by Theorem \ref{ZD27Aug20}. Hence, the algebra $A/(P,Q)$ is also not  a commutative domain.

 It remains to consider the case when the pair  of ideals is $(x_i)$, $P$ where $P\in \Spec_{x_1x_2}(A)$. Suppose that the algebra $A$ is of type (C10). Then, by Theorem \ref{18Sep20}.(1), $P=\gp A$ for some prime ideal $\gp \in \Spec Z(A))\backslash V(T_A)$. 
  By Theorem \ref{D27Aug20},  the algebra $$A/(x_i, P)=\bigg(\bigoplus_{i\in I, \beta \in \CR}\Big(Z(A)/\gp\Big)x_3^ix^\beta\bigg)/(x_i)$$ is not commutative. 
  
 Suppose that the algebra $A$ is of type (C11) and the algebra $A/(x_i,P)$ is a commutative domain. 
 We seek a contradiction. The factor algebras
 \begin{eqnarray*}
A/(x_1)&=&K\langle x_2,x_3\, | \, x_2x_3=(x_3-\mu )x_2\rangle \simeq 
\mA ,\\
A/(x_2)&=&K\langle x_1,x_3\, | \, x_1x_3=(x_3-\alpha )x_1\rangle \simeq 
\mA, 
\end{eqnarray*}
are both isomorphic to the skew polynomial ring $\mA = K[x_3][x_1, \; \s ]$ where $\s (x_3)=x_3-1$. By the classification of prime ideals of the algebra $\mA$, every prime ideal $\gq$ of the algebra $\mA$ such that $\mA/ \gq$ is a commutative algebra contains the element $x_1$.  Therefore, $x_2\in (x_1,P)$ (resp., $x_1\in (x_2,P$)). In a view of the $(x_1,x_2)$-symmetry, let us consider the first case. Recall that $A_{x_1x_2}=T^{-1}_AA$ (Lemma \ref{XA15Sep20})  where $T_A=M_\S$ (Lemma \ref{a15Sep20}). Then, 
by Theorem \ref{18Sep20}.(2), 
$$x_2=x_1a+(x_1^nx_2^n)^{-s}(bx_1^{pn}+\sum_{i=1}^{p-1}b_ix_1^{\xi_in}x_2^{in}+cx_2^{pn})$$
for some elements $a,b,c, b_i\in A$ and a natural number $s$.  By multiplying the equality above on the left by the element 
 $(x_1^nx_2^n)^s$ and using the $\N^2$-grading of the algebra $A$ we get a contradiction. $\Box$ 

\begin{theorem}\label{A27Aug20}
Let $A=A(q, \alpha ,\mu)$ and $C:=C_A(x_1,x_2)$, see Proposition \ref{D24Sep19}. Then 
\begin{enumerate}
\item Suppose that $q^2\neq 1$. Then 
$$\Aut_K(A)=\{ \s_\l \,| \, \l= (\l_1,\l_2,\l_3)\in K^\times \times K^\times \times C\}\simeq K^\times \times K^\times \times C, \;\; \s_\l\mapsto \l$$ where
$\s_\l : x_1\mapsto \l_1x_1$, $x_2\mapsto \l_2x_2$, $x_3\mapsto x_3+\l_3$.
\item Suppose that $q^2= 1$ ($q\neq 1$). Then 
$$\Aut_K(A)=\{ \s_\l \tau^i\,| \, \l= (\l_1,\l_2,\l_3)\in K^\times \times K^\times \times C, i=0,1\}\simeq \Big(K^\times \times K^\times \times C\Big)\rtimes \Z_2,$$ a semidirect product of groups where 
$\s_\l : x_1\mapsto \l_1x_1$, $x_2\mapsto \l_2x_2$, $x_3\mapsto x_3+\l_3$ and $\tau : x_1\mapsto x_2$, $x_2\mapsto x_1$, $x_3\mapsto x_3$; $\tau^2=1$, $\tau \s_{(\l_1,\l_2,\l_3)}=\s_{(\l_2,\l_1,\tau (\l_3))}\tau$ and $\Z_2=\Z/2 \Z\simeq \langle \tau \rangle$. 
\end{enumerate}

\end{theorem}

{\it Proof}. Suppose that $q^2\neq 1$ (resp., $q^2=1$).  Let $\G =\{ \s_\l \,| \, \l= (\l_1,\l_2,\l_3)\in K^\times \times K^\times \times C\}$. Using the defining relations of the algebra $A$, we see that $\G$ (resp., $\G':= \{ \s_\l \tau^i\,| \, \l= (\l_1,\l_2,\l_3)\in K^\times \times K^\times \times C, i=0,1\}$) is a subgroup of $\Aut_K(A)$ which is isomorphic  to the group $ K^\times \times K^\times \times C$ via $\s_\l\mapsto \l$ (resp., $\G \rtimes \Z_2$). 

Given an automorphism $\s$ of the algebra $A$. We have to show that $\s \in \G$ (resp., $\s \in \G'$). The theorem follows from the following Claim.

{\sc Claim:} {\em The group $\Aut_K(A)$ acts on the set $\{ K^\times x_1, K^\times x_2\}$ by the rule $(\s , S)\mapsto \s (S)$ where $S\in \{ K^\times x_1, K^\times x_2\}$}.

Suppose that the Claim holds. Then we have two options either $\s (x_1)=\l_1x_1$, $\s (x_2)=\l_2x_2$ or $\s (x_1)=\l_1x_2$, $\s (x_2)=\l_2x_1$ for some elements $\l_1, \l_2\in K^\times$. The second option is only possible in the case when $q^2=1$. This follows from the relation $x_2x_1=qx_1x_2$.
 In the second case, by replacing the automorphism $\s$ by $\s\tau$ we may assume that in both cases $\s (x_1)=\l_1x_1$ and  $\s (x_2)=\l_2x_2$.  Furthermore, by replacing $\s$ by the automorphism $\s \s^{-1}_{(\l_1, \l_2, 0)}$, we may assume that $\l_1=\l_2=1$. Let $x_3'=\s (x_3)$. Then 
 \begin{eqnarray*}
 [x_3'-x_3,x_1]&=&\s ([x_3,x_1])-[x_3,x_1]=\s  (\alpha x_1)-\alpha x_1=\alpha x_1-\alpha x_1=0,  \\
 \, [x_3'-x_3,x_2] &=& \s ([x_3,x_2])-[x_3,x_2]=\s  (\mu x_2)-\mu x_2=\mu x_2-\mu x_2=0,
\end{eqnarray*}
 and so $c:=x_3'-x_3\in C$, i.e. $\s =\s_{(1,1,c)}\in \G$, as required.
 
{\it Proof of the Claim}. Suppose that the algebra $A$ is of type (C1)--(C7) (resp.,  (C10) or (C11)). 
By Proposition \ref{a21Sep20}.(1) (resp., Proposition \ref{a21Sep20}.(3)), the group $\Aut_K(A)$ acts on the set $\{ (x_1), (x_2)\}$. Since  $K^\times$ is the group of units of the algebra $A$  the Claim follows.

Suppose that the algebra $A$ is of type (C8) or (C9). In the view of $(x_1,x_2)$-symmetry, it suffices to consider only the case (C8), i.e.  $q$ is  a primitive $n$'th  root of unity, $p>0$, $\alpha =1$, $\mu =0$ and $p\nmid n$.  Clearly,  the algebra $A=K[x_2,x_3][x_1; \s]$ is a skew polynomial ring  where $\s(x_2)=q^{-1}x_2$ and $\s (x_3)=x_3-1$.  By Proposition \ref{a21Sep20}.(2), the ideal $(x_1)$ is the only height 1 prime ideal $P$ of the algebra $A$ such that $A/P$ is a commutative domain. Hence, $\s ((x_1))=(x_1)$, and so $\s (x_1) =\l_1 x_1$ for some unit $\l_1 $ of the algebra $A$, i.e. $\l_1 \in K^\times$.

   By Theorem \ref{D27Aug20}, $A=\bigoplus_{i\in I, \beta \in \CR} Z(A) x_3^ix^\beta$ where $I=\N_{<p}$ and $\CR = \N_{<pn}\times \N_{<n}$. Notice that 
 $$ \bigoplus_{\beta_1\in N_{<pn},\beta_2\in \N_{<n}}Z(A)x_1^{\beta_1}x_2^{\beta_2}$$ is the direct sum of the eigenspaces of the automorphism $\o_{x_1}$ where the set of eigenvalues is $\{ q^{-\beta_2}\, | \, \beta_2\in N_{<n}\}$. In particular, $$\s^{\pm 1}\Bigg(\bigoplus_{\beta_1\in N_{<pn},\beta_2\in \N_{<n}}Z(A)x_1^{\beta_1}x_2\Bigg)=\bigoplus_{\beta_1\in N_{<pn},\beta_2\in \N_{<n}}Z(A)x_1^{\beta_1}x_2.$$ Hence, $\s (x_2)=\l_2x_2$ for some element $\l_2\in K^\times$, and the Claim holds.  $\Box $\\

{\bf Classification of normal and $\o$-weight elements of the algebra $A$.}  
 Theorem \ref{20Sep20} is a classification of $\o$-weight elements of the algebra $A$. It is also a classification of normal elements of the algebra $A$ as we will show that every normal element of $A$ is an $\o$-weight element and vice versa (Proposition \ref{A20Sep20}).

\begin{theorem}\label{20Sep20}

\begin{enumerate}
\item Suppose that the algebra $A$ belongs to one of the cases (C1)--(C6), (C7b) or (C8)--(C10). Then $\CW_A(\o )=\{ (q^{-\beta_2}, q^{\beta_1}, \alpha\beta_1+\mu \beta_2)\, | \, \beta \in \CR_A\}$ and for each weight $\l = (q^{-\beta_2}, q^{\beta_1}, \alpha\beta_1+\mu \beta_2)$, $A_\l = Z(A)x^\beta$ where the set $\CR_A$ is defined in Theorem \ref{D27Aug20}.
\item Suppose that the algebra $A$ belongs to the case  (C7a). Then  $\CW_A(\o )=\{ (q^{-\beta_2}, q^{\beta_1}, \beta_1+\mu \beta_2+kn)\, | \, \beta \in \N^2_{<n}, k\in \N\}$ and for each weight $\l = (q^{-\beta_2}, q^{\beta_1}, \beta_1+\mu \beta_2+kn)\in \CW_A(\o )$, $$A_\l = \bigoplus_{i=0}^{k}Kx_1^{\beta_1 +in}x_2^{\beta_2+(k-i)n}\;\; {\rm  and }\;\; \dim_K(A_\l )=k+1<\infty.$$ In particular, the weight spaces are finite dimensional. 
\item Suppose that the algebra $A$ belongs to the case  (C11). Then  $\CW_A(\o )=\CW_{A_{x_1x_2}}(\o )=\{ (q^{-\beta_2}, q^{\beta_1}, \beta_1+\mu \beta_2)\, | \, \beta \in \CR_{x_1x_2}=\N_{<pn}\times \N_{<n} \}$ and for each weight $\l_\beta  = (q^{-\beta_2}, q^{\beta_1}, \beta_1+\mu \beta_2)\in \CW_A(\o )$, 
$$A_{\l_\beta} = \bigoplus_{i=0}^{p-1}\L'x^{\th (\beta , i)}$$ 
where $\L'= K[x_3^p-x_3][x_1^{pn}, x_2^{ pn}]$ and $\th (\beta , i)$ is a unique element of the set $\N_{<pn}^2$ such that $\th (\beta , i)\equiv \beta +in(-\mu , 1)\mod pn$.

\end{enumerate}
\end{theorem}

{\it Proof}. 1. By Corollary \ref{a28Aug20}.(1), $\Es_A(\o ) =\Es_A(\o_{x_1},\o_{x_2})= \bigoplus_{\beta \in \CR_A}Z(A)x^\beta$, and  for all $z\in Z(A)$ and $\beta \in \CR_A$, $\o_{x_1}(zx^\beta )=q^{-\beta_2}zx^\beta$,  $\o_{x_2}(zx^\beta )=q^{\beta_1}zx^\beta$, and $\ad_{x_3}(zx^\beta)=(\alpha\beta_1+\mu \beta_2)zx^\beta$. Statement 1 follows from the Claim.

{\sc Claim.} {\em The map $\CR_A\ra \CW_A(\o )$, $\beta \mapsto \l_\beta =(q^{-\beta_2}, q^{\beta_1}, \alpha \beta_1+\mu \beta_2)$ is an injection.}

Suppose that the Claim holds. Then necessarily the map is  well defined and a bijection since $\Es_A(\o ) = \bigoplus_{\beta \in \CR_A}Z(A)x^\beta$, and this direct sum is the direct sum of $\o$-eigenspaces. 

{\it Proof of the Claim}. 

Suppose that the algebra $A$ is of type (C1)--(C3). The the Claim is true since the element $q$ is not a root of unity.

Suppose that the algebra $A$ is of type (C4) or (C5). In view of Theorem \ref{A28Oct18}.(1b), it suffices to consider the case (C4). Recall that in this case, $q$ is a primitive $n$'th  root of unity, $p=0$, $\alpha =1$, $\mu =0$ and $\CR_A=\N\times \N_{n}$ (Theorem \ref{D27Aug20}). Then $\l_\beta =(q^{-\beta_2}, q^{\beta_1}, \beta_1)= \l_{\beta'}=(q^{-\beta_2'}, q^{\beta_1'}, \beta_1')$ implies that $\beta_2=\beta_2'$ and  $\beta_1=\beta_1'$, and so $\l_\beta  = \l_{\beta'}$. 

Suppose that the algebra $A$ is of type (C6). In this case, $q$ is a primitive $n$'th  root of unity, $p=0$, $\alpha =1$ and $\mu =-\frac{\mu_1}{\mu_2}$ where $\mu_1,\mu_2\in \N_+$ such that $\mu_1\neq \mu_2$ and $(\mu_1,\mu_2)=1$, and $\CR_A= \mS (\mu_1n, \mu_2n)=\N^2\backslash \Big( (\mu_1n, \mu_2n)+\N^2\Big)$ (Theorem \ref{D27Aug20}). Then 
$\l_\beta =(q^{-\beta_2}, q^{\beta_1}, \beta_1+\mu \beta_2)= \l_{\beta'}=(q^{-\beta_2'}, q^{\beta_1'}, \beta_1'+\mu \beta_2')$ implies the equality $\mu_2 (\beta_1'-\beta_1)=\mu_1(\beta_2'-\beta_2)$, and so $\beta'=\beta +k(\mu_1, \mu_2)$ for some $k\in \N$ since $(\mu_1, \mu_2)=1$, i.e. $\beta' =\beta $ (since $\beta,\beta'\in \mS (\mu_1n, \mu_2n) $), as required. 

Suppose that the algebra $A$ is of type (C7b). In this case,  $\alpha =1$, $\mu\not\in \Q $ and $\CR_A=\N^2$ (Theorem \ref{D27Aug20}). Then $\l_\beta =(q^{-\beta_2}, q^{\beta_1}, \beta_1+\mu \beta_2)= \l_{\beta'}=(q^{-\beta_2'}, q^{\beta_1'}, \beta_1'+\mu \beta_2')$ implies $\beta = \beta'$ since $\mu \not\in \Q$ (use the equality of the last coordinates). 

Suppose that the algebra $A$ is of type (C8) or (C9). In view of the obvious symmetry, it suffices to consider the type (C8). 
In this case, $q$ is  a primitive $n$'th  root of unity, $p>0$, $\alpha =1$, $\mu =0$, $p\nmid n$ and $\CR_A=\N_{<pn}\times \N_{<n}$  (Theorem \ref{D27Aug20}). Then $\l_\beta =(q^{-\beta_2}, q^{\beta_1}, \beta_1)= \l_{\beta'}=(q^{-\beta_2'}, q^{\beta_1'}, \beta_1')$ implies $\beta = \beta'$ by comparing the first and the last coordinates of the vectors $\l_\beta$ and $\l_{\beta'}$.

Suppose that the algebra $A$ is of type (C10). In this case, $q$ is  a primitive $n$'th  root of unity, $p>0$, $\alpha =1$, $\mu\not\in \mF_p $, $p\nmid n$ and $\CR_A=\N_{<pn}\times \N_{<pn}$ (Theorem \ref{D27Aug20}). Then $\l_\beta =(q^{-\beta_2}, q^{\beta_1}, \beta_1+\mu \beta_2)= \l_{\beta'}=(q^{-\beta_2'}, q^{\beta_1'}, \beta_1'+\mu \beta_2')$ implies 
$\beta = \beta'$ since $\mu \not\in \Fp$ (use the equality of the last coordinates). 

2. By Corollary \ref{a28Aug20}.(1), $\Es_A(\o ) =\Es_A(\o_{x_1},\o_{x_2})= \bigoplus_{\beta \in \CR_A}Z(A)x^\beta$, and  for all $z\in Z(A)$ and $\beta \in \CR_A$, $\o_{x_1}(zx^\beta )=q^{-\beta_2}zx^\beta$,  $\o_{x_2}(zx^\beta )=q^{\beta_1}zx^\beta$, and $\ad_{x_3}(zx^\beta)=(\alpha\beta_1+\mu \beta_2)zx^\beta$. Recall that in the case (C7a), $q$ is a primitive $n$'th  root of unity, $p=0$, $\alpha =1$,  $\mu =\frac{\mu_1}{\mu_2}$,
 where $\mu_1,\mu_2\in \N_+$ such that  $(\mu_1,\mu_2)=1$, $Z(A)=K$ (Theorem \ref{C24Sep19})   and $\CR_A=\N^2$ (Theorem \ref{D27Aug20}). Then $\l_\beta =(q^{-\beta_2}, q^{\beta_1}, \beta_1+\mu \beta_2)= \l_{\beta'}=(q^{-\beta_2'}, q^{\beta_1'}, \beta_1'+\mu \beta_2')$ iff $\beta_i=\g_i +l_in$,   $\beta_i'=\g_i +l_i'n$ for unique $\g_i\in N_{<p}$ and  for some $l_i, l_i'\in N$ where  $i=1,2$,    
 and $\g_1+\mu \g_2 +(l_1+l_2)n=\g_1+\mu \g_2 +(l_1'+l_2')n$. It follows that 
 $\CW_A(\o )=\{ (q^{-\g_2}, q^{\g_1}, \g_1+\mu \g_2+kn)\, | \, (\g_1\g_2) \in \N^2_{<n}, k\in \N\}$ and for each weight $\l = (q^{-\g_2}, q^{\g_1}, \g_1+\mu \g_2+kn)\in \CW_A(\o )$, $$A_\l = \bigoplus_{i=0}^{k}Kx_1^{\g_1 +in}x_2^{\g_2+(k-i)n}\;\; {\rm  and }\;\; \dim_K(A_\l )=k+1<\infty.$$

3.  By Corollary  \ref{b16Sep20}, $\CW_{A_{x_1x_2}}(\o )=\{ (q^{-\beta_2}, q^{\beta_1}, \beta_1+\mu \beta_2)\, | \, \beta \in \CR_{x_1x_2}=\N_{<pn}\times \N_{<n} \}$ and 
for each weight $\l_\beta  = (q^{-\beta_2}, q^{\beta_1}, \beta_1+\mu \beta_2)\in \CW_{A_{x_1x_2}}(\o )$, $A_{x_1x_2, \l_\beta}=Z(A_{x_1x_2})x^\beta$. Notice that $\Es_A(\o ) \subseteq A\cap \Es_{A_{x_1x_2}}(\o )$. Hence, 
  to prove statement 3, it suffices to show that 
$$A\cap A_{x_1x_2,\l_\beta} = \bigoplus_{i=0}^{p-1}\L'x^{\th (\beta , i)}$$ 
(since then $\CW_A(\o )=\CW_{A_{x_1x_2}}(\o )$ and for each element  $\l_\beta\in \CW_A(\o )$, $A_{\l_\beta} =A\cap A_{x_1x_2,\l_\beta}$).
\begin{eqnarray*}
 A\cap A_{x_1x_2,\l_\beta} &=& A\cap Z(A_{x_1x_2})x^\beta =A\cap K[x_3^p-x_3][x_1^{\pm pn}, (x_1^{-\mu n}x_2^n)^{\pm 1},x_2^{\pm pn}]x^\beta\\
  &=&A\cap\bigoplus_{i=0}^{p-1}K[x_3^p-x_3][x_1^{\pm pn}, x_2^{\pm pn}](x_1^{-\mu n}x_2^n)^ix^\beta\\
   &=& A\cap\bigoplus_{i=0}^{p-1}K[x_3^p-x_3][x_1^{\pm pn}, x_2^{\pm pn}]x^{\beta+in(-\mu , i)}\\
 &=& \bigoplus_{i=0}^{p-1}\L'x^{\th (\beta , i)}
\end{eqnarray*}
where $\th (\beta , i)$ is a unique element of the set $\N_{<pn}^2$ such that $\th (\beta , i)\equiv \beta +in(-\mu , 1)\mod pn$. $\Box $

\begin{proposition}\label{A20Sep20}
Every normal element of the algebra $A$ is an $\o$-weight vector, and vice versa.
\end{proposition}

{\it Proof}. $(\Rightarrow )$ Let $a=\sum_{i\in \N, \beta \in \N^2}\l_{i\beta}x_3^ix^\beta$ be a nonzero normal element of the algebra $A$ where $\l_{i\beta}\in K$. Then $ab=\o_a(b) b$ for all elements $b\in A$ where $\o_a$ is the automorphism of the algebra $A$ that os determined by the normal nd regular element $A$. In particular, for $i=1,2$, 
$$\o_a(x_i)a=ax_i=x_i\o_{x_i}^{-1}(a),$$
and so $\o_{x_i}^{-1}(a)=c_ia$ where $c_i=x_i^{-1}\o_a(x_i)\in A_{x_i}$. Let $(\alpha_1, \alpha_2):=(\alpha , \mu )$. On the other hand, $\o_{x_i}^{-1}(a)=\sum_{i\in \N, \beta \in \N^2}\l_{i\beta}(x_3+\alpha_i)^i\rho_{i,\beta}x^\beta$ where 
$$\rho_{i,\beta}:=\begin{cases}
q^{\beta_2}& \text{if }i=1,\\
q^{-\beta_1}& \text{if }i=2.\\
\end{cases}
$$
Using the $\N^2$-grading of the algebra $A$ we see that $c_i\in K^\times$. 

Similarly, $ \o_a(x_3)a=ax_3=\sum_{i\in \N, \beta \in \N^2}(x_3-(\alpha\beta_1+\mu\beta_2))\l_{i\beta}x_3^ix^\beta$. Using the $\N^2$-grading of the algebra $A$ we see that
$\alpha\beta_1+\mu\beta_2$ is a constant, say $\nu$,  for all $i\in \N$ and $ \beta \in \N^2$ such that $\l_{i\beta}\neq 0$. So, 
$\ad_{x_3}(a)=\nu a$. Therefore, the element $a$ is an $\o$-weight vector.

$(\Leftarrow )$ Every  $\o$-weight vector of the algebra $A$ is a normal element: If $a$ is an $\o$-eigenvector with eigenvalue $\l_\beta$ then 
$x_1a=q^{-\beta_2}ax_1$, $x_2a=q^{-\beta_1}ax_2$ and $ax_3=(x_3-\alpha \beta_1-\mu \beta_2)a$. $\Box $\\



\section{The Krull, global and Gelfand-Kirillov dimensions of the algebras $\CA = \bigotimes_{i=1}^n A_i$ where $A_i=A(q_i, \alpha_i ,\mu_i)$}\label{KRULLGLDIM}

The aim of this section is to compute the Krull,   classical Krull,  global and Gelfand-Kirillov dimensions of the algebras  $\CA = \bigotimes_{i=1}^n A_i$ where $A_i=A(q_i, \alpha_i ,\mu_i)$. The centre of the algebra $\CA$ is described. \\

{\bf The algebras $\CA = \bigotimes_{i=1}^n A_i$ where $A_i=A(q_i, \alpha_i ,\mu_i)$.}

\begin{proposition}\label{A23Aug20}
Let $\CA = \bigotimes_{i=1}^n A_i$ be a tensor product of algebras $A_i=A(q_i, \alpha_i ,\mu_i)=K\langle x_{i1}, x_{i2}, x_{i3}\rangle$  in  Theorem \ref{A28Oct18}.(1).

\begin{enumerate}
\item The algebra $\CA $ is an iterated skew polynomial algebra, 
$$\CA = K[x_{11}][x_{12}; \s_{12}][x_{13}; \s_{13}]\cdots [x_{n1}][x_{n2}; \s_{n2}][x_{n3}; \s_{n3}].$$

\item The algebra $\CA $ is a Noetherian domain with PBW basis, $\CA = \bigoplus_{\beta \in \N^{3n}}Kx^\beta$.   
\item The  Gelfand-Kirillov dimension of the algebra $\CA$ is $3n$.
\item $Z(\CA ) = \bigotimes_{i=1}^n Z(A_i)$.
 
\end{enumerate}
\end{proposition}

{\it Proof.} 1. Statement 1 follows at once from Theorem \ref{A28Oct18}.(1,a).

2. By statement 1, the algebra $\CA$ is a Noetherian domain with PBW basis (since $\CA = \otimes_{i=1}^n A_i$).

3. Statement 3 follows at once from statement 2. 

4. This statement is obvious (as the centre of a tensor product of algebras is the tensor product of their centres). $\Box $\\

{\bf $\N^{2n}$-grading of the algebra $\CA =\bigotimes_{i=1}^nA_i$.} By (\ref{AN2gr}),
the algebra $\CA$
is an $\N^{2n}$-graded algebra,
\begin{equation}\label{AN2gr2}
\CA =\bigoplus_{\beta \in \N^{2n}} P_nx^\beta \;\; {\rm where}\;\; P_n=K[x_{13}, \ldots , x_{n3}]\;\; {\rm and}\;\; x^\beta = x_{11}^{\beta_1}x_{12}^{\beta_2}\cdots x_{n1}^{\beta_{2n-1}}x_{n2}^{\beta_{2n}}.
\end{equation}
The $\N^{2n}$-grading on the algebra $\CA$ is the tensor product of $\N^2$-grading of the tensor multiples $A_i$ of the algebra $\CA$. The algebra $P_n$ is a polynomial algebra in $n$ variables. 

The product in the algebra $A$ is given by the rule: For all elements $a,b\in P_n$ and $\beta , \g \in \N^{2n}$,
$$ ax^\beta \cdot bx^\g= a \o^\beta (b) q^{[\beta ,\g]} x^{\beta +\g}\;\; {\rm where}\;\; \o^\beta =\prod_{i=1}^n\o_{x_{i1}}^{\beta_{2i-1}}\o_{x_{i2}}^{\beta_{2i}}\;\; {\rm and}\;\; q^{[\beta ,\g]}=\prod_{i=1}^n q_i^{\beta_{2i}\g_{2i-1}} .$$
So, the algebra $\CA$ is a skew crossed product of the (additive) monoid $\N^{2n}$ with  the polynomial algebra $P_n$.  \\

 {\bf The global, Krull and classical Krull dimensions of algebras $\CA$.} For an algebra $A$, we denote by $\CK (A)$ and  cl.Kdim $(A )$  the Krull and the classical Krull dimension of $A$, respectively.

\begin{theorem}\label{A24Sep19}
Let $\CA = \bigotimes_{i=1}^n A_i$ be a tensor product of algebras $A_i=A(q_i, \alpha_i ,\mu_i)=K\langle x_{i1}, x_{i2}, x_{i3}\rangle$  in  Theorem \ref{A28Oct18}.(1).

\begin{enumerate}
\item The global dimensional of $\CA$ is $3n$.
\item The Krull dimensional of $\CA$ is $3n$.
\item The classical Krull dimensional of $\CA$ is $3n$.
 \end{enumerate}
\end{theorem}

{\it Proof.} 1. By Proposition \ref{A23Aug20}.(1), the algebra $\CA$ is an iterated skew polynomial algebra. Now, the result follows from \cite[Theorem 7.5.3.(iii)]{MR}. 

2. By Proposition \ref{A23Aug20}.(1), the algebra $\CA$ is an iterated skew polynomial algebra. Now, the result follows from \cite[Proposition 6.5.4.(i)]{MR}.

3. (i)  cl.Kdim $(\CA )\geq 3n$: By Theorem \ref{A28Oct18}.(1),  there is a strictly ascending chain of primes in the algebra $\CA$: 
$$ 0\subset \gp_1\subset \cdots \subset \gp_{3n}=(x_{11}, x_{12}, x_{13}, \ldots , x_{n1}, x_{n2}, x_{n3})$$
where for each natural number $k=1, \ldots , n-1$,
$$\gp_{3k+i}=\begin{cases}
(x_{11}, x_{12}, x_{13}, \ldots , x_{k1}, x_{k2}, x_{k3})& \text{if }i=0,\\
(x_{11}, x_{12}, x_{13}, \ldots , x_{k1}, x_{k2}, x_{k3}, x_{k+1, 1})& \text{if }i=1,\\
(x_{11}, x_{12}, x_{13}, \ldots , x_{k1}, x_{k2}, x_{k3}, x_{k+1, 1},x_{k+1, 2} )& \text{if }i=2,\\
\end{cases}
$$
and the statement (i) follows. 

(ii)  cl.Kdim $(\CA )\leq 3n$: The algebra $\CA$ is a Noetherian algebra (Proposition \ref{A23Aug20}.(2)), hence  cl.Kdim $(\CA )\leq \CK (A)=3n$, by \cite[Lemma 6.4.5]{MR}.

By the statements (i) and (ii), cl.Kdim $(\CA )= 3n$.  $\Box$

{\bf Licence.} For the purpose of open access, the author has applied a Creative Commons Attribution (CC BY) licence to any Author Accepted Manuscript version arising from this submission.\\

{\bf Disclosure statement.} No potential conflict of interest was reported by the author.\\

{\bf Data availability statement.} Data sharing not applicable – no new data generated.

\small{

\begin{tabular}{l  l}

V. V.   Bavula  \quad \quad \quad \quad  \quad \quad \quad \quad \quad \quad \quad  & Ali Al Khabyah  \\ 
Department of Pure Mathematics \quad \quad \quad \quad  \quad \quad \quad \quad \quad \quad \quad  & Department of Mathematics\\
University of Sheffield \quad \quad \quad \quad  \quad \quad \quad \quad \quad \quad \quad  & Faculty of Science \\
Hicks Building \quad \quad \quad \quad  \quad \quad \quad \quad \quad \quad \quad  &  Jazan University \\
Sheffield S3 7RH \quad \quad \quad \quad  \quad \quad \quad \quad \quad \quad \quad  & Jazan 45142 \\
UK \quad \quad \quad \quad  \quad \quad \quad \quad \quad \quad \quad  &  Saudi Arabia \\
email: v.bavula@sheffield.ac.uk \quad \quad \quad \quad  \quad \quad \quad \quad \quad \quad \quad  & email: aalkhabyah@jazanu.edu.sa\\
\end{tabular}

}

\end{document}